\documentclass[11pt, epsfig]{article}
\usepackage{epsfig, amsmath, amssymb, amsthm, times}
\usepackage{graphicx,color}

\newcommand{\R}{{\mathbb R}}
\newcommand{\C}{{\mathbb C}}
\newcommand{\N}{{\mathbb N}}
\newcommand{\Z}{{\mathbb Z}}

\def\bW{{\mathbf W}}

\DeclareMathOperator{\Cay}{Cay}

\usepackage[mathcal]{eucal}

\newtheorem {thm}{Theorem}[section]

\theoremstyle{defintion}
\newtheorem {df}[thm]{Definition}

\theoremstyle{remark}

\theoremstyle{example}
\newtheorem{ex}[thm]{Example}

\theoremstyle{assumption}

\def\P{\operatorname{\mathbb P}}
\def\T{{\mathbb T}}

\def\W{\mathcal{W}}
\def\Gn{\Gamma_n}
\def\lbl{\label}
\def\be{\begin{equation}}
\def\ee{\end{equation}}
\def\p{\partial}

\newcommand{\1}{{i\mkern1mu}}

\title{Bifurcations in the Kuramoto model on graphs}
\author{Hayato Chiba,\thanks{Institute of Mathematics for Industry, 
Kyushu University / JST PRESTO, Fukuoka,
819-0395, Japan, {\tt chiba@imi.kyushu-u.ac.jp}}\;
Georgi S. Medvedev,\thanks{Department of Mathematics, 
Drexel University, 3141 Chestnut Street, Philadelphia, PA 19104,
{\tt medvedev@drexel.edu}} \;
and Matthew S. Mizuhara\thanks{Department of 
Mathematics and Statistics,
The College of New Jersey,
{\tt  mizuharm@tcnj.edu}}
}
\begin{document}
\maketitle

\begin{abstract}
In his classical work, Kuramoto analytically
described the onset of synchronization in all-to-all coupled networks of phase oscillators with
random intrinsic frequencies. Specifically, he identified a critical value of the coupling strength,
at which the incoherent state loses stability and a gradual build-up of coherence begins.
Recently, Kuramoto's scenario was shown to hold for a large class of coupled systems on convergent 
families of deterministic and random graphs \cite{ChiMed17a, ChiMed17b}. 
Guided by these results, in the present work, we study several model problems
illustrating the link between network topology and synchronization in coupled dynamical systems.

First, we identify several families of graphs, for which the transition to synchronization in 
the Kuramoto model starts at the same critical value of the coupling strength and proceeds in practically
the same way. These examples include Erd\H{o}s-R\' enyi random graphs, Paley graphs, complete bipartite 
graphs, and certain stochastic block graphs. 
These examples illustrate that some rather simple structural properties such as the volume of the graph 
may determine the onset of synchronization, while finer structural features may affect only higher order statistics
of the transition to synchronization. Further, we study the transition to synchronization in the Kuramoto model on 
power law and small-world random graphs. The former family of graphs endows the Kuramoto model with 
very good synchronizability: the synchronization threshold can be made arbitrarily low by varying the parameter of the power
law degree distribution. For the Kuramoto model on small-world graphs, in addition to the transition to synchronization,
we identify a new bifurcation leading to stable random twisted states. The examples analyzed in this work complement 
the results in \cite{ChiMed17a, ChiMed17b}.
\end{abstract}

\section{Introduction}
\setcounter{equation}{0}

The Kuramoto model of coupled phase oscillators provides an important paradigm for studying collective dynamics and 
synchronization in ensembles of interacting 
dynamical systems. In its original form, the KM describes the dynamics of all-to-all coupled 
phase oscillators with randomly distributed intrinsic frequencies
\be\lbl{cKM}
\dot \theta_{i} = \omega_{i} + Kn^{-1} \sum_{j=1}^n  \sin(\theta_{j}-\theta_{i}),
\quad i\in [n]:=\{1,2,\dots,n\}.
\ee
Here, $\theta_i:\R\to \T:=\R/2\pi\Z$ is the phase of the oscillator $i\in [n]$. The intrinsic
frequencies $\omega_i,\; i\in [n],$ are independent identically distributed random variables.
The density of the probability distribution of $\omega_1$, $g,$ is a smooth even function.
The sum on the right-hand side of \eqref{cKM} models the interactions between the oscillator
$i$ and the rest of the network. Parameter $K$ controls the strength of the interactions.
The sign of $K$ determines the type of interactions. The coupling is attractive 
if $K>0$, and is repulsive otherwise. Sufficiently strong attractive coupling favors synchrony.

\begin{figure}
\begin{center}
\textbf{a} \includegraphics[height=1.6in,width=1.8in]{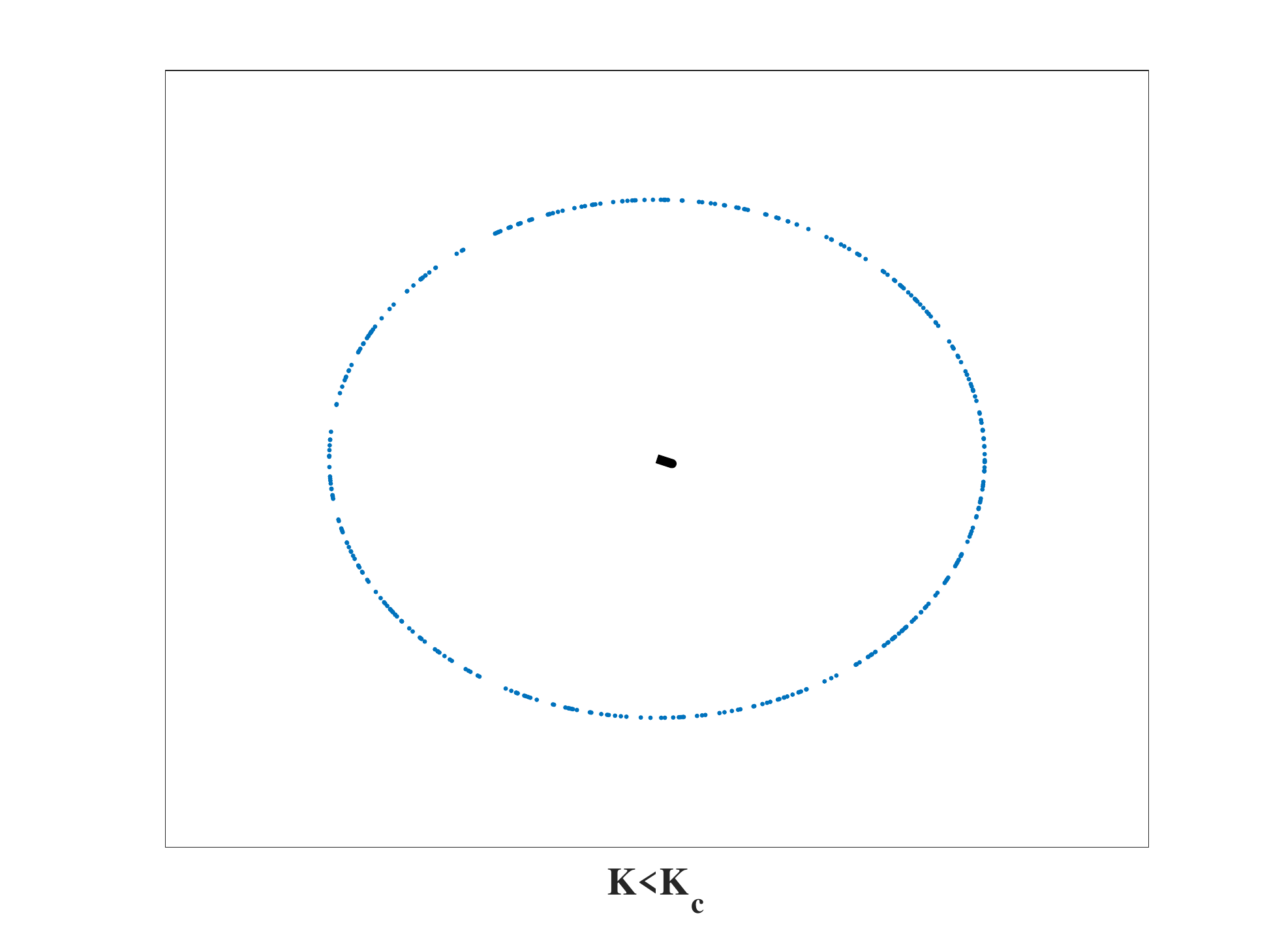} 
\textbf{b} \includegraphics[height=1.6in,width=1.8in]{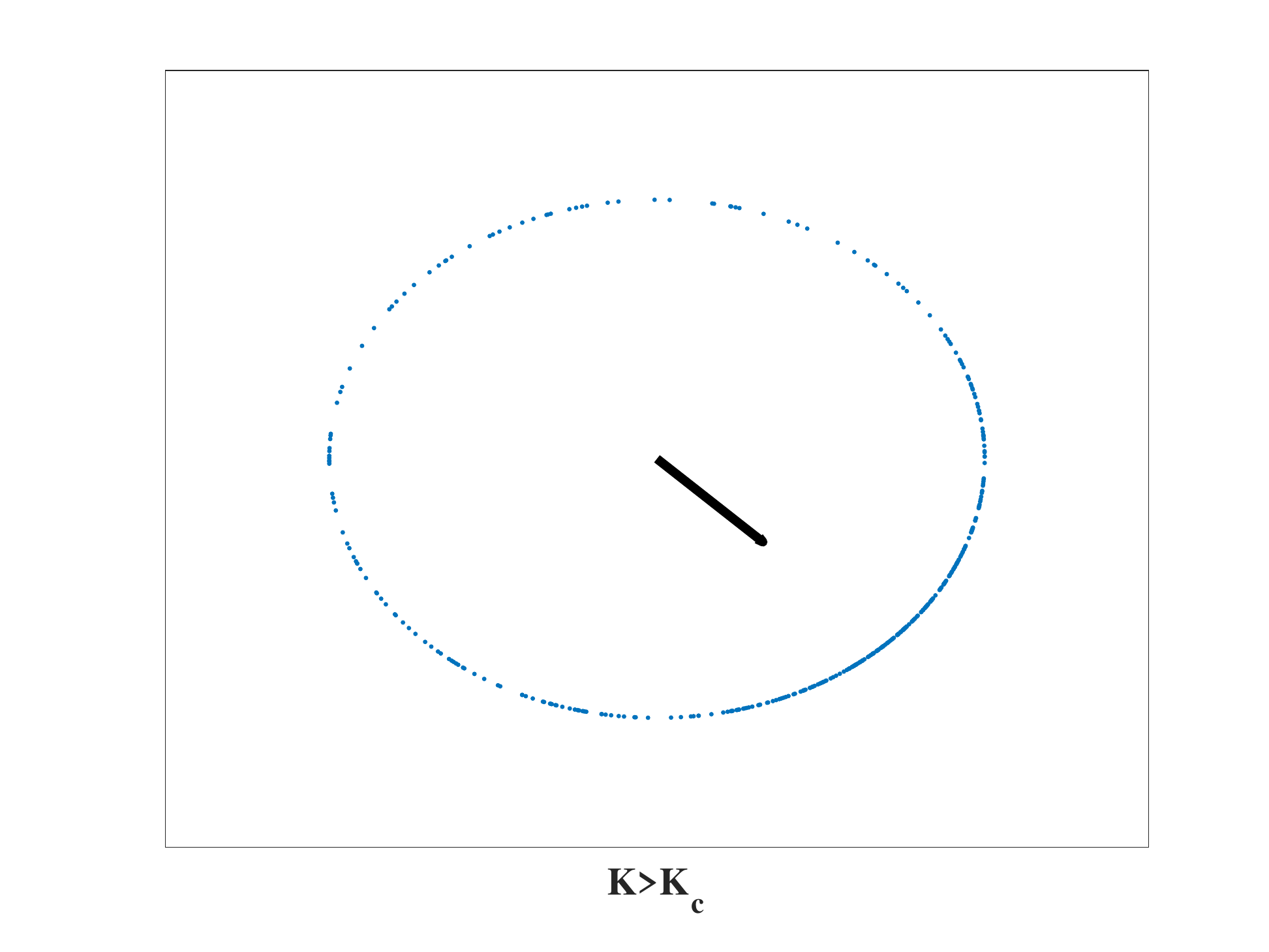}
\textbf{c} \includegraphics[height=1.6in,width=1.8in]{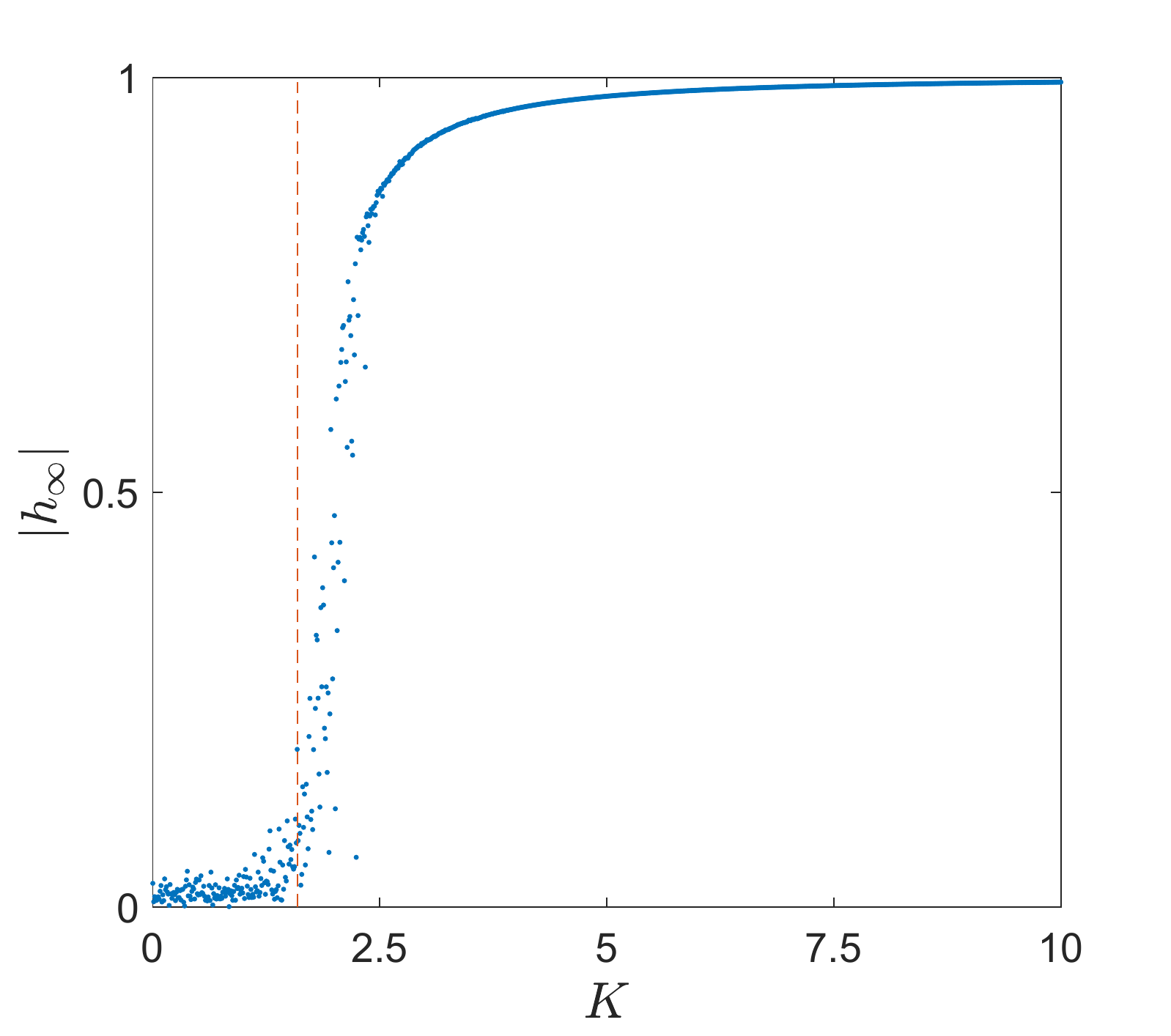}
\caption{The distribution of the phases of coupled oscillators is shown on the unit circle 
in the complex plane: $\theta_{k}\mapsto e^{\1 \theta_k}\in \C,\; k\in [n]$. The strength of 
coupling is below the critical value $K_c$ in (\textbf{a}) and  is above $K_c$ in  (\textbf{b}).  
The black arrow depicts the order parameter, as a vector in the complex plane (cf.~\eqref{order}).
The bigger size (modulus) of the order parameter corresponds to the higher degree of coherence.
In (\textbf{a}) the modulus of the order parameter is close to zero, and the distribution of the 
oscillators is close to the uniform distribution. In contrast in (\textbf{b}), the distribution develops
a region of higher density. The complex order parameter points to the
expected value (center of mass) of the distribution of the oscillators around $\T$.
 \textbf{c}) The modulus of the order parameter is plotted for
 different values of $K$. The gradual change of the modulus of the
 order parameter from values around zero to those close to $1$
marks the transition to synchronization.
}
\label{f.1}
\end{center}
\end{figure}

For  small positive values of $K$, the KM shows little coherence. The phases fill out the unit circle 
approximately uniformly (Fig.~\ref{f.1}\textbf{a}). This dynamical regime is called an incoherent state. 
It persists for positive
values of $K$ smaller than the critical value $K_c=2\left(\pi g(0)\right)^{-1}$. For values of $K>K_c$, the
system undergoes a gradual build-up of coherence approaching complete synchronization as 
$K\to\infty$ (Fig.~\ref{f.1}\textbf{b}). Kuramoto identified the critical value $K_c$ and described the 
transition to synchrony, using the complex \textit{order parameter} 
\be\lbl{order}
h(t)=n^{-1}\sum_{j=1}^n e^{\1 \theta_j(t)}.
\ee
Below, we will often use the polar form of the order parameter
\be\lbl{polar-order}
h(t)=r(t)e^{\1\psi (t)}.
\ee

Specifically, he showed that for $t\gg 1$  (cf.~\cite{Str00})
\be\lbl{transition}
r(t)= \left\{ \begin{array}{ll} O(n^{-1/2}),& 0<K<K_c,\\
r_\infty (K-K_c)+O(n^{-1/2}), & K>K_c,
\end{array}
\right.
\ee
where
\be\lbl{r-inf}
r_\infty(x)=r_1\sqrt{x} + O(x),\; x\ge 0, \qquad r_1={4\over K_c^2\sqrt{ -\pi g^{\prime\prime}(0)}}.
\ee
Here, $0\le r(t)\le 1$ and $\psi(t)$ stand for the modulus and the argument of the order parameter defined
by the right-hand side of \eqref{order}.  The value of $r$ is interpreted as a measure of coherence in 
the system dynamics. Indeed, if all phases are independent random variables distributed uniformly on $\T$ 
(complete incoherence) then  $r=o(1)$ with probability one, by the Strong Law of Large Numbers.
If, on the other hand, all phase variables are equal then $r=1$. Equations \eqref{transition}, \eqref{r-inf} 
suggest that the 
system undergoes a pitchfork bifurcation en route to synchronization.  This bifurcation is clearly seen 
in the numerical experiments (see Figure~\ref{f.1}\textbf{c}). Furthermore, the Kuramoto's scenario of the transition
to synchrony was recently confirmed by rigorous mathematical analysis of the KM \cite{Chi15a,Die16,FGG16}.

\begin{figure}
\begin{center}
\textbf{a} \includegraphics[height=1.6in,width=2.0in]{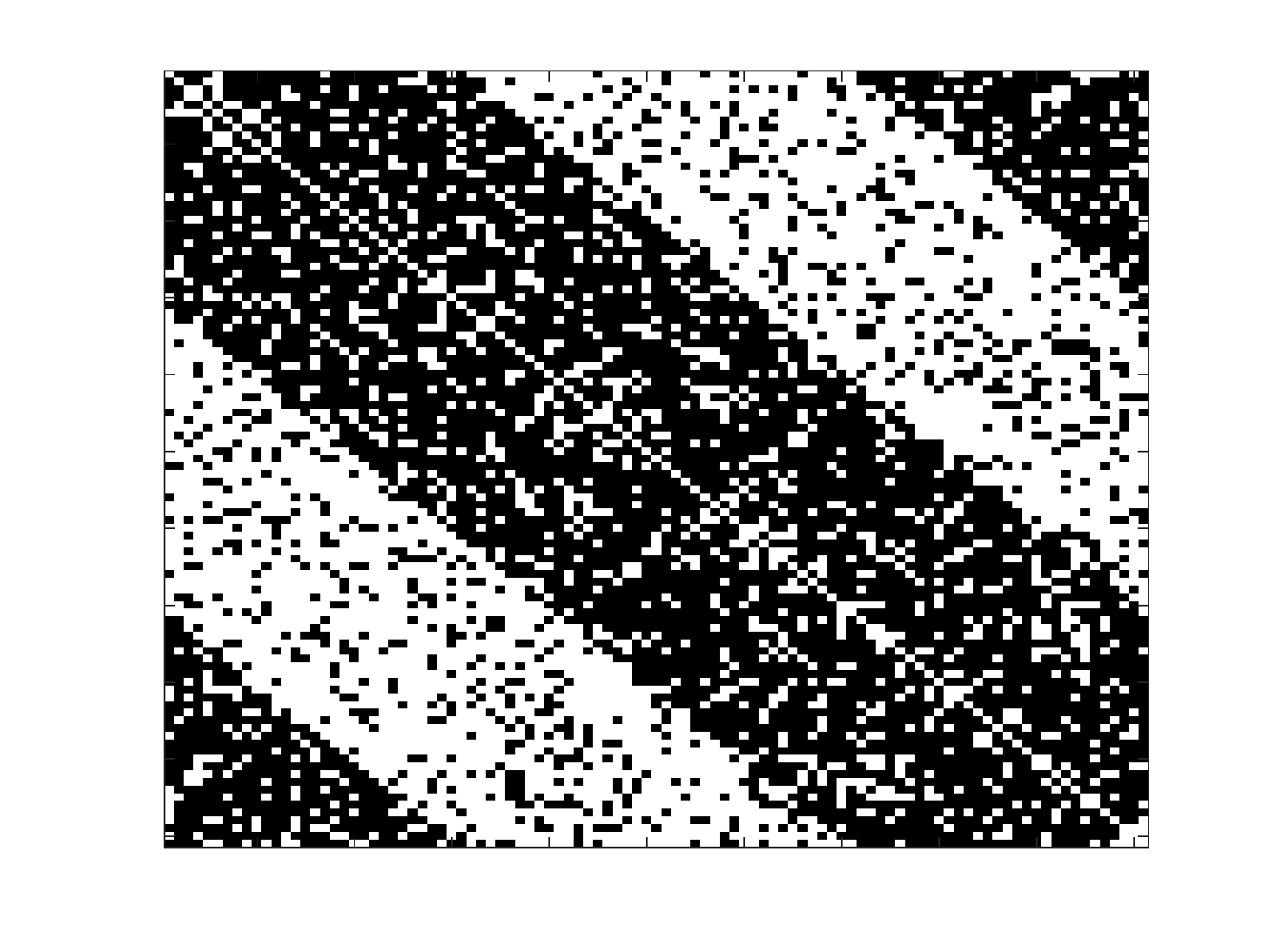}
\textbf{b} \includegraphics[height=1.6in,width=2.0in]{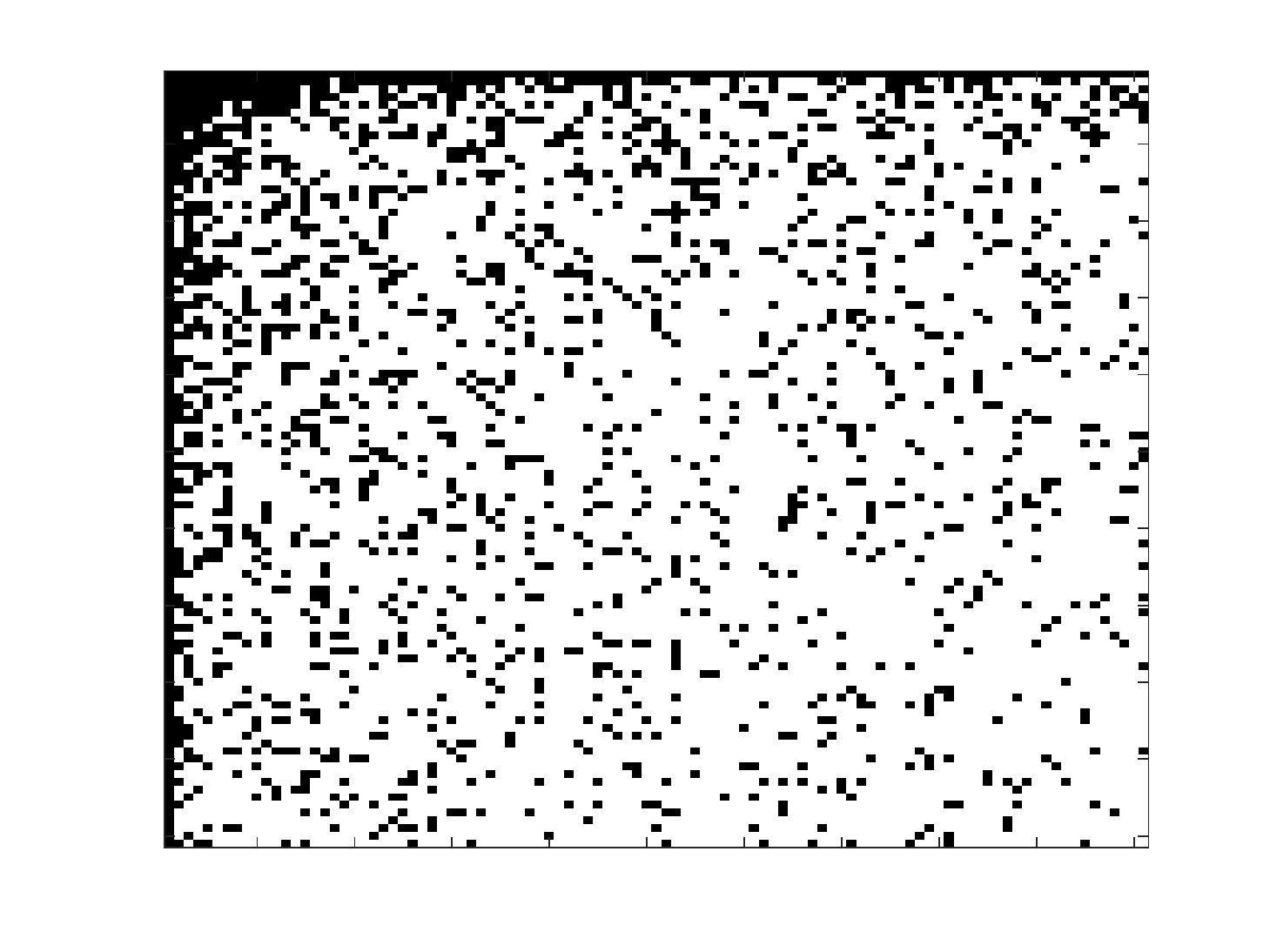}
\caption{ Pixel pictures of small-world (\textbf{a}) and power law (\textbf{b})
graphs on $100$ nodes.
}
\label{f.graphs}
\end{center}
\end{figure}

In this paper, we investigate  bifurcations in the KM on a variety of graphs ranging from symmetric Caley
graphs to random small-world and power law graphs. To visualize the connectivity of a large graph, we 
use the pixel picture of a graph. This is a square black and white plot, where each black pixel stands 
for entry $1$ in the adjacency matrix of the graph (see Figure~\ref{f.graphs}). As we  show below, the graph 
structure plays a role in the transition to synchrony. Furthermore, some graphs may force bifurcations
to spatial patterns other than synchrony.
To highlight the role of the network topology, we present numerical
experiments of the KM on small-world graphs \cite{WatStr98}. These graphs are formed by 
replacing some of the
edges of a regular graph with random edges (see Figure~\ref{f.graphs}\textbf{a}). 
These results demonstrate that some graphs may exhibit bifurcations to spatial patterns other 
than synchrony.

\begin{figure}
\begin{center}
\textbf{a} \includegraphics[height=1.6in,width=2.0in]{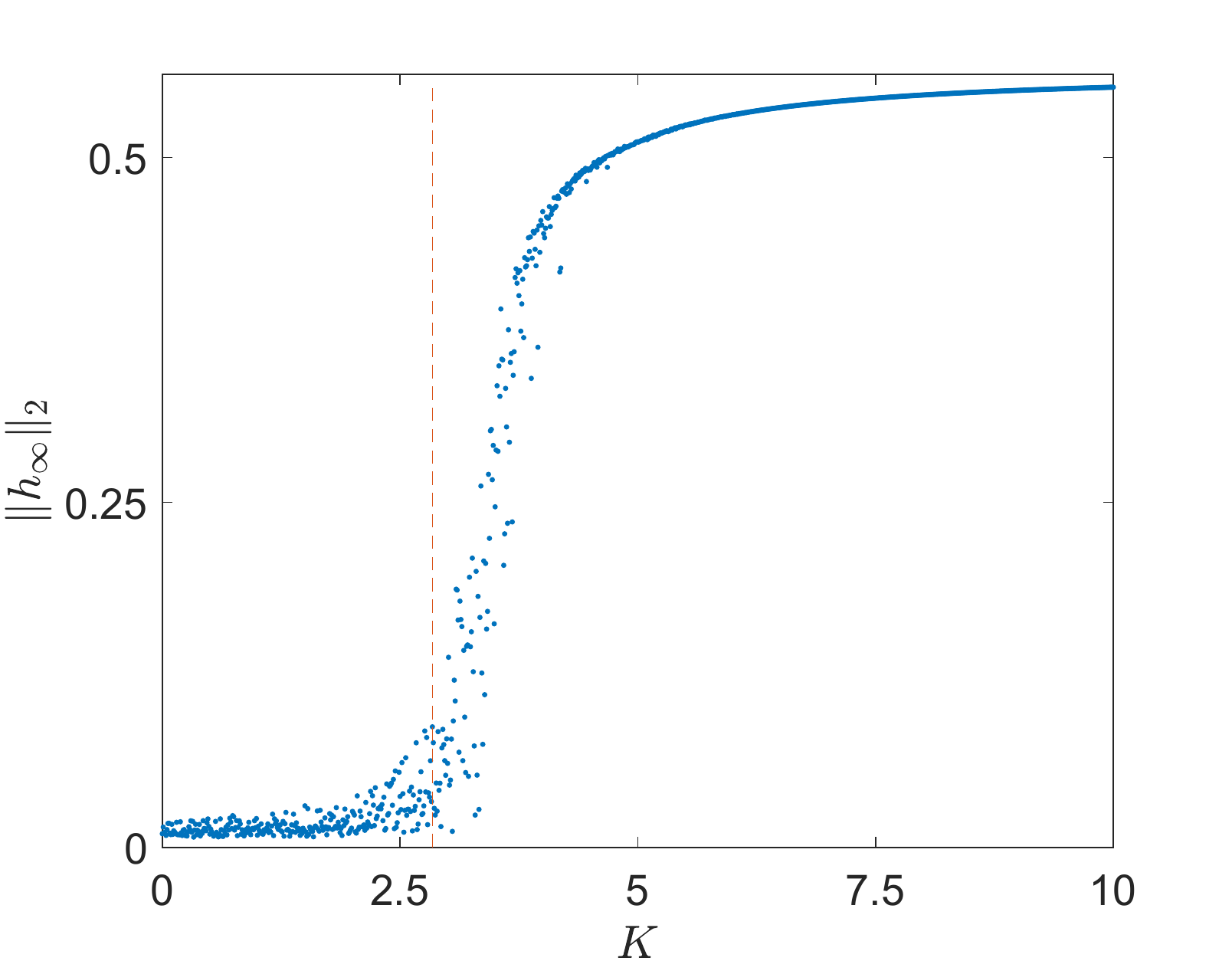}
\textbf{b} \includegraphics[height=1.6in,width=2.0in]{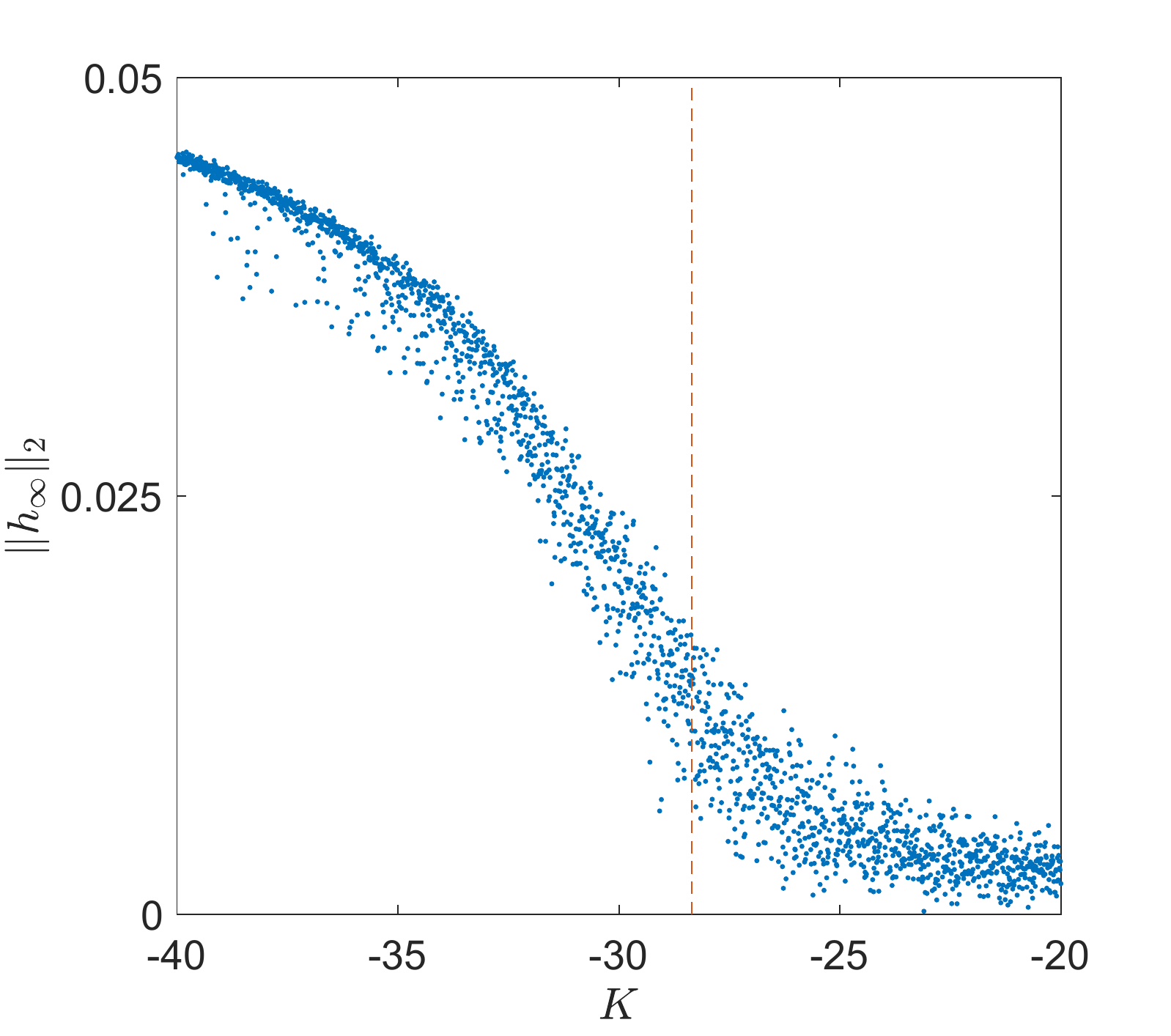}\\
\caption{Transition to synchronization in the KM on small-world graphs (\textbf{a}).
For negative $K$, the KM  undergoes another transition at $K_c^-<0$ (\textbf{b}). 
The emerging pattern is shown in Figure~\ref{f.2+} \textbf{a},\textbf{b}. 
}
\label{f.2}
\end{center}
\end{figure}

\begin{figure}
\begin{center}
\textbf{a} \includegraphics[height=1.6in,width=2.0in]{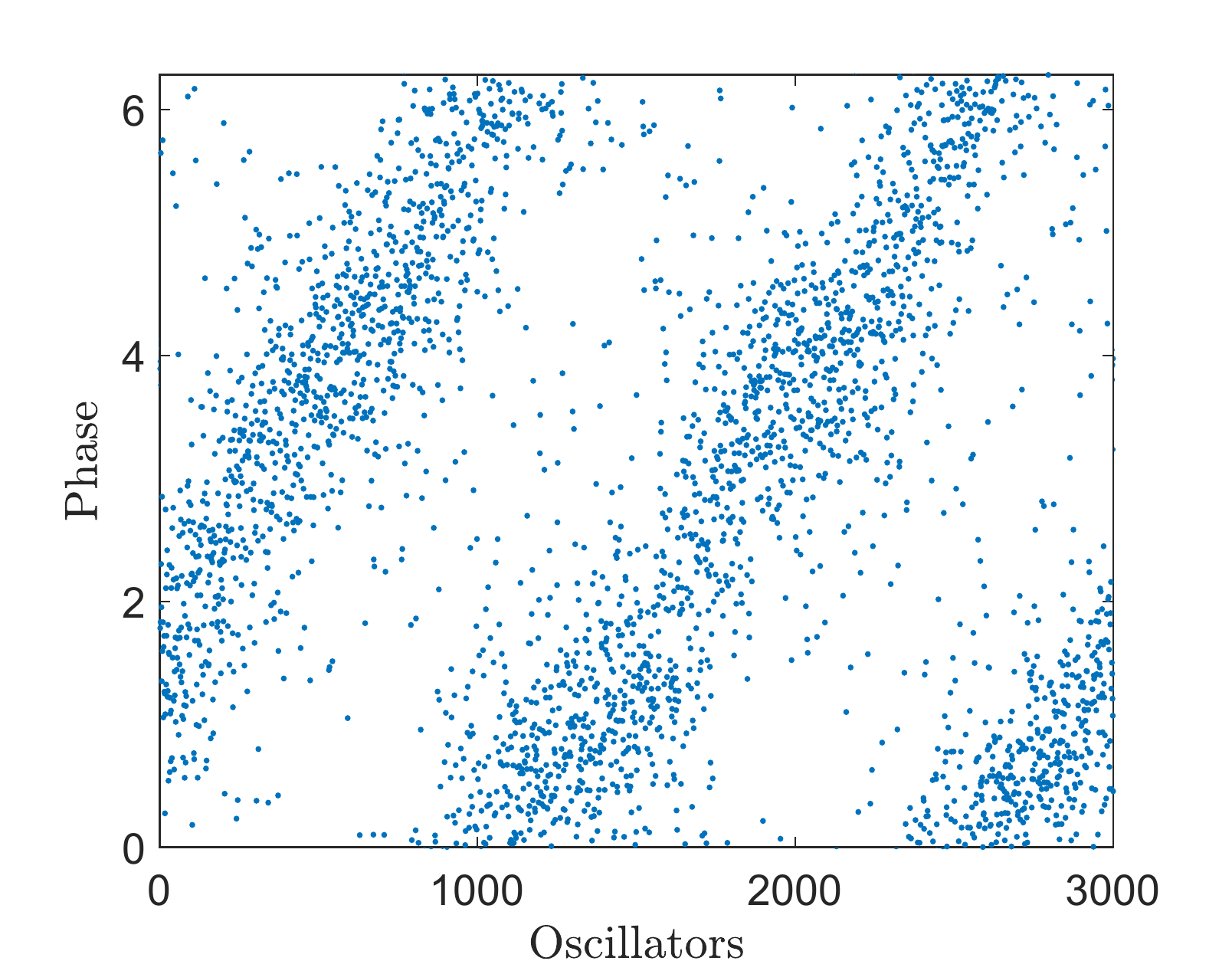}
\textbf{b} \includegraphics[height=1.6in,width=2.0in]{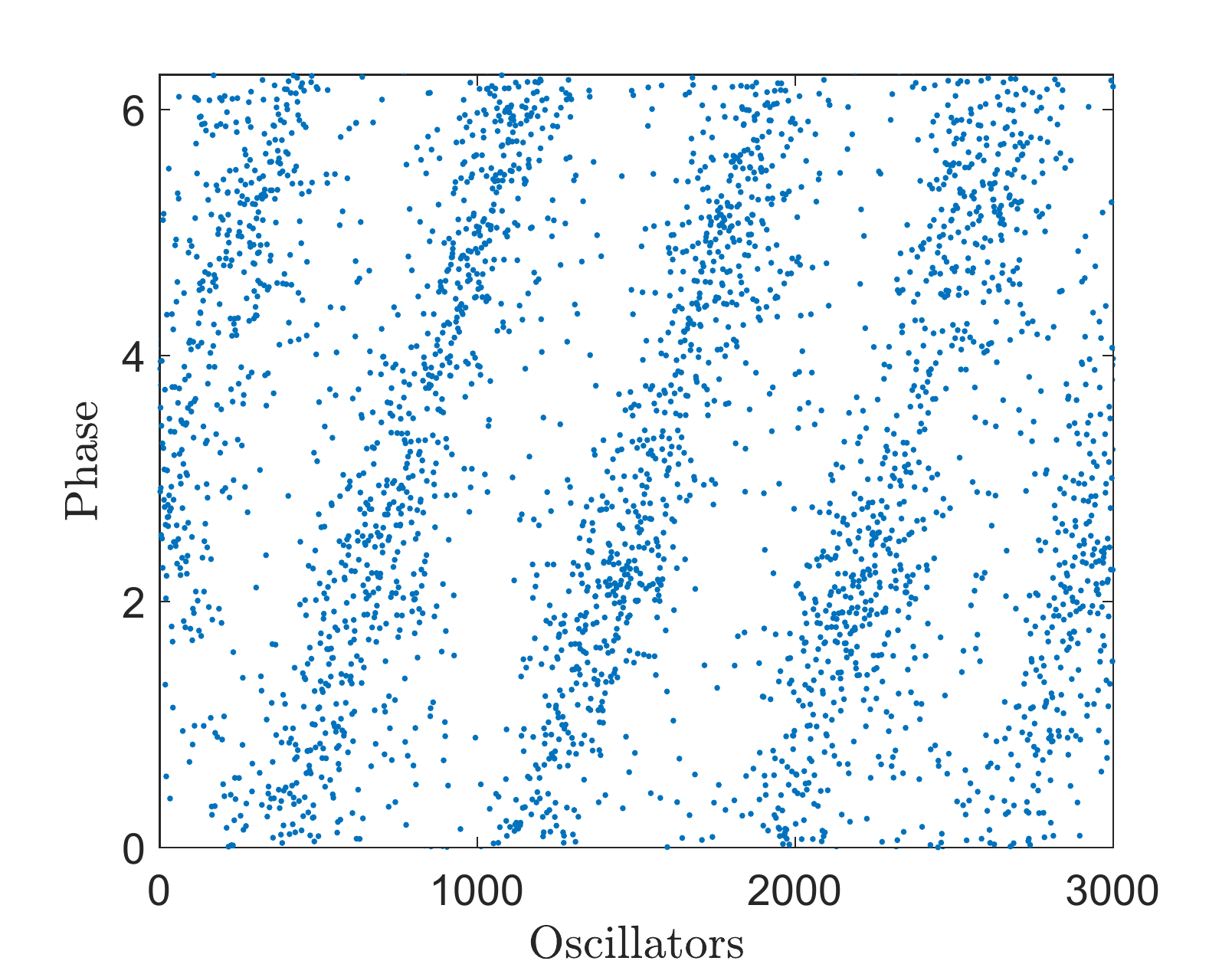}
\caption{The oscillators in \textbf{a} are 
clustered around a $2$-twisted state ($p=.2, r=.3, K=-36$; see
Example~\ref{ex.Wrandom} for the definitions of parameters). By changing one of the parameters controlling small-world
connectivity (namely, by reducing the range of local interactions), 
one can obtain $q$-twisted states for different values of $q\in\Z$, such 
as the $4$-twisted state shown in plot \textbf{b} ($p=.2, r=.2, K=-40$).
}
\label{f.2+}
\end{center}
\end{figure}

Figure~\ref{f.2}\textbf{a} shows that for the KM on small-world graphs like in that on complete graphs, 
the order parameter  undergoes a smooth transition\footnote{For the KM on graphs, the order 
parameter must be suitably redefined (see \eqref{dWorder}).}. Near 
the critical value $K_c^+\approx 3.2$, the asymptotic in time value of the order parameter starts to grow monotonically 
and approaches $0.5$ for increasing values of $K$. This change in the order parameter corresponds to the 
transition to synchronization. Interestingly, in contrast to the original KM, the model on small-world graphs
exhibits another transition at a negative value of $K$, $K_c\approx -27$ (Figure~\ref{f.2}\textbf{b}).   This time it grows 
monotonically for decreasing values of $K$ and approaches the value approximately equal to $0.0475$. The 
corresponding state of the network is shown in Figure~\ref{f.2+}\textbf{a}. In the pattern shown in this figure, 
the oscillators are distributed randomly about a $2$-twisted state. 
A $q$-twisted state ($q\in\Z$) is a linear function on the unit circle:
$T_q(x)= {2\pi qx\over n} \mod 2\pi,$ $x\in [n],$ $q$-twisted states have been studied before as stable steady states in repulsively
coupled KM on Caley graphs \cite{MedTan15b, WilStr06}. By changing parameters controlling small-world 
connectivity, we find
different $q$-twisted states bifurcating from the incoherent state for decreasing negative $K$ 
(see Figure~\ref{f.2+}\textbf{a,b}). These numerical experiments show that in the KM on small-world graphs
the incoherent state is stable for $K\in (K_c^-, K_c^+)$. For values $K>K_c^+,$ the system undergoes the 
transition to synchrony, while for $K<K^-_c$, it leads to random patterns localized around $q$-twisted states.

The loss of stability of the incoherent state in the KM on spatially structured networks was analyzed in 
\cite{ChiMed17a, ChiMed17b}. The linear stability analysis of the incoherent state yields the critical 
values $K_c^\pm$ \cite{ChiMed17a}, and the analysis of the bifurcations at $K^\pm_c$ explains
 the emerging spatial patterns \cite{ChiMed17b}. In the present work,  using the insights from
\cite{ChiMed17a, ChiMed17b}, we conduct numerical 
experiments illustrating the role of the network topology in shaping the transition from the 
incoherent state to coherent structures like synchronous and twisted states in the KM on graphs. 
In particular, we present numerical simulations showing for many distinct graphs the onset of synchronization
takes place at the same critical value and the graph structure has only higher order effects. We identify
the dominant structural properties of the graphs shaping the transition to synchronization in the KM.

The effective analysis of the KM on large graphs requires an analytically convenient description of 
graph sequences.
To this end, we employ the approach developed in \cite{Med14a, Med14b} inspired by the theory of graph limits
\cite{LovGraphLim12} and, in particular, by W-random graphs \cite{LovSze06}. We explain the graph models 
used in this paper in Section~\ref{sec.KM}.
The W-random graph framework affords an analytically tractable mean field description of the KM on a broad
class of dense and sparse graphs. The mean field partial differential equation for the KM on graphs  is explained
in Section~\ref{sec.MF}.
There we also explain the generalization of the order parameter suitable for the KM on graphs. After  that we 
go over the main results of \cite{ChiMed17a, ChiMed17b}. In Section~\ref{sec.results}, we present our 
numerical results. We use carefully designed examples of graphs to highlight the structural properties 
affecting the onset of synchrony. For instance, we show that the KMs on Erd\H{o}s-Renyi and Paley graphs
have the same critical value $K_c$ marking the onset of synchronization. Furthermore, the mean field limits
of the KM on these graphs coincide with that for the KM on weighted complete graph.
We demonstrate the similarity in the transition to synchronization in the KM on bipartite graphs and on the 
family of stochastic block graphs interpolating a disconnected graph and a weighted complete graph.
This time the mean field equations for the two models are different, but the critical value $K_c$ remains 
the same. Further, we investigate the bifurcations in the KM on small-world and power law graphs illuminating
the effects of these random connectivity patterns on the dynamics of the KM.

\section{The KM on graphs}\lbl{sec.KM}
\setcounter{equation}{0}
We begin with the description of graph models that will be used in this paper.
Let $\Gn=\langle V(\Gn), E(\Gn), A_n\rangle$ be a weighted directed graph on $n$ nodes.
$V(\Gn)=[n]$ stands for the node set of $\Gn$.  $A_n=(a_{n,ij})$ is an
$n\times n$ weight matrix. The edge set
$$
E(\Gn)=\left\{ \{i,j\}\in [n]^2:\;a_{n,ij}\ne 0\right\}.
$$
If $\Gn$ is a simple graph then $A_n$ is the adjacency matrix.

The KM on $\Gn$ is defined as follows
\begin{eqnarray}
\lbl{KM}
\dot u_{n,i} &=& \omega_i + K(n\alpha_n)^{-1} \sum_{j=1}^n a_{n,ij} \sin (u_{n,j}-u_{n,i}),
\quad i\in [n],\\
\lbl{KM-ic}
u_{n,i}(0)&=&u_{n,i}^0. 
\end{eqnarray}
The sequence $\alpha_n\ne 0$ is only needed if $\{\Gn\}$ is a sequence of sparse graphs.
Without proper rescaling, the continuum limit
(as $n\to\infty$) of the KM on sparse graphs is trivial. The sequence $\{\alpha_n\}$  will be explained below, when
we introduce sparse $W$-random graphs. Until then one can can assume $\alpha_n\equiv 1$.

We will be studying solutions of \eqref{KM} in the limit as $n\to\infty$. Clearly, such limiting behavior is 
possible only if the sequence of graphs $\{\Gn\}$ is convergent in the appropriate sense. We want to define
$\{\Gn\}$ in such a way that it includes graphs of interest in applications, and at the same time is convenient
for deriving the continuum limit. To achieve this goal, we use the ideas from the theory of graph
limits \cite{LovGraphLim12}. Specifically, we choose a square integrable function $W$ on a unit
square $I^2:=[0,1]^2$. $W$ is used to define the asymptotic behavior of $\{\Gn\}$ as $n\to\infty$.
It is called a graphon in the theory of graph limits \cite{LovGraphLim12}.
Let 
$$
\W:=\left\{ W\in L^2(I^2):\quad W(x,y)=W(y,x)   \;\mbox{a.e.}\right\}.
$$
We define $\W_0\subset \W$ as a set of measurable symmetric functions on $I^2$ with values in $I$.
We discretize the unit interval $I$: $x_{n,j}=j/n, \; j\in \{0\}\cup [n]$ and denote
$
I_{n,i}:=(x_{n,i-1}, x_{n,i}], \; i\in [n].
$

The following graph models
are inspired by W-random graphs 
\cite{LovSze06, BCCZ}. Let $\Gn=\langle V(\Gn), E (\Gn), A_n=(a_{n,ij})\rangle$. 
\begin{description}
\item[(DD)] Weighted deterministic graph $\Gn= H(n,W)$ is defined as follows: 
\be\lbl{a-det}
a_{n,ij} =\langle W\rangle_{I_{n,i}\times I_{n,j}}:= n^2 \int_{I_{n,i}\times I_{n,j}} W(x,y)dxdy,
\ee
where $W\in\W$.
\item[(RD)] Dense random graph $\Gn= H_r(n,W)$. Let  $W\in\W_0$ and 
$a_{n,ij},\; 1\le i\le j\le n,$ be independent identically distributed (IID) random variables such that 
\be\lbl{Pdense}
\P\left(  a_{n,ij}=1 \right)=\langle W\rangle_{I_{n,i}\times I_{n,j}},\quad 
\P\left(  a_{n,ij}=0 \right)=1-\P\left(  a_{n,ij}=1 \right).
\ee 
\item[(RS)] Sparse random graphs $\Gn= H_r(n,W, \{\alpha_n\})$. Let $W\in\W$
be a nonnegative function and positive sequence $1\ge \alpha_n\searrow
0$  satisfy $n\alpha_n\to \infty$ as $n\to\infty$, 
\be\lbl{Psparse}
\P\left( \{j,i\} \in E(\Gamma_n) \right)=\alpha_n \langle \tilde W_n \rangle_{I_{n,i}\times I_{n,j}}, \quad 
\tilde W_n(x,y):= \min\{\alpha_n^{-1}, W(x,y)\}.
\ee
\end{description}

We illustrate the random graph models in \textbf{(RD)} and \textbf{(RS)} with the following examples.
\begin{ex}\lbl{ex.Wrandom}
\begin{description}
\item[(SW)] Small-world graph $S_{n,r,p}=H_r(n,W_{r,p})$, $r\in (0,0.5)$, $p\in (0,0.5],$ is defined via
\be\lbl{def-Wpr}
W_{p,r}(x,y)=
\left\{ 
\begin{array}{ll}
1-p, & \min \{|x-y|, 1-|x-y|\} \le r,\\
p, &\mbox{otherwise},
\end{array}
\right.
\ee
(see Figure~\ref{f.graphs}\textbf{a}).
\item[(ER)] Erd\H{o}s-R{\' e}nyi graph $G_{n,p}=H_r(n, W_p)$,
$W_p\equiv p\in (0,1)$.
\item[(SER)] Sparse  Erd\H{o}s-R{\' e}nyi graph $\tilde G_{n,p}=H_r(n, W_p, \{ n^{-\beta}\})$, \; $\beta\in (0,1)$.
\item[{(PL)}] Power law graph $\Gamma= H_r(n, W^{\gamma}, \{n^{-\beta}\})$ 
\be \lbl{Walpha}
W ^\gamma (x,y)= (xy)^{-\gamma}, \quad 0<\gamma<\beta<1,
\ee
(see Figure~\ref{f.graphs}\textbf{b}).
\end{description}
\end{ex}

Graphon $W$ carries all information needed for the derivation of the continuum limit for the KM on
deterministic and random graphs (\textbf{DD}), (\textbf{RD}), and (\textbf{RS}). However,
for the KM on random graphs \textbf{(RD)} and \textbf{(RS)} taking the continuum limit involves 
an additional step. Namely, we first approximate the model on random graph by that on 
the averaged deterministic weighted 
graph $\Gamma_n=H(n,W)$:
\begin{eqnarray}\lbl{aKM}
\dot v_{n,i} & =&  \omega_{n,i}+Kn^{-1}\sum_{j=1}^n \bar W_{n,ij} \sin(v_{n,j}-v_{n,i}),\quad
i\in [n],\\
v_{n,i}(0) &=& u^0_{n,i},
\end{eqnarray}
where 
$
\bar W_{n,ij}:=\langle \tilde W_n\rangle_{I_{n,i}\times I_{n,j}}.
$

The approximate model \eqref{KM} is formally derived from the original KM \eqref{KM} on a random graph
$\Gn=H_r(n,W, \{ \alpha_n\})$ by averaging the right hand side of \eqref{KM} over all possible realizations 
of $\Gn$. The justification of averaging is given in \cite{KVMed17a}.

\section{The mean field limit}\lbl{sec.MF}
\setcounter{equation}{0}

Solutions of the KM with distributed frequencies such as incoherent
state and the solutions bifurcating from the incoherent state are best described in statistical terms.
To this end, suppose 
$\rho(t,u,\omega,x)$ is the conditional density of the random vector $(u,\omega)$
 given $\omega$, and parametrized by $(t,x)\in \R^+\times I$.
Here,  the spatial domain $I=[0,1]$ represents the continuum of the oscillators in the limit $n\to\infty$
(see \cite{ChiMed17a} for precise meaning of the continuum limit). Then $\rho$ satisfies the 
following initial value problem
\begin{eqnarray}
\lbl{MF}
{\p \over \p t} \rho(t,u,\omega,x) +
{\p \over \p u }\left\{\rho(t,u,\omega,x) V(t,u,\omega,x)\right\}
&=&0,\\
\lbl{MF-ic}
\rho(0,u,\omega, x) &=& \rho^0(u),
\end{eqnarray}
where initial density $\rho^0(u)$ independent of $\omega$ and $x$ for
simplicity (see \cite{KVMed17a} for the treatment of a more general case). 
The velocity field $V$ is defined via
\be\lbl{def-V}
V(t,u,\omega,x)=\omega+
K\int_I\int_\R\int_\T W(x,y) \sin(\phi - u ) \rho(t,\phi, \omega ,y) g(\omega ) 
d\phi d\omega d y.
\ee

The density $\rho(t,\cdot)$ approximates the distribution of the oscillators around $\T$ at time $t\in [0,T]$.
Specifically, in the large $n$ limit the empirical measures defined on the Borel $\sigma$-algebra
$\mathcal{B}(G),$ $G:=\T\times \R\times I$:
\be\lbl{EMeasure}
\mu_t^n(A)=n^{-1} \sum_{i=1}^n \delta_{(u_{ni}(t), \omega_i, x_{ni})} (A), \quad  A\in \mathcal{B}(G),
\ee
converge weakly to the absolutely continuous measure
\be\lbl{CMeasure}
\mu_t(A)=\int_A \rho(t,u,\omega, x)g(\omega ) du d\omega dx, 
\ee
provided the initial data \eqref{KM-ic} converge weakly to $\mu_0$ (cf.~\cite[Theorem~2.2]{ChiMed17a}, see also
\cite{KVMed17a}).

The mean field limit provides a powerful tool for studying the KM. It gives a simple analytically
tractable description of complex dynamics in this model. In particular, the incoherent state in the mean field
description corresponds to the stationary density $\rho_u=1/2\pi$. The 
stability analysis of $\rho_u,$ a steady state solution of \eqref{MF}, yields the region of stability of the
 incoherent state and the critical values of $K$, at which it loses stability. Furthermore, the analysis of the 
bifurcations at the critical values of $K$ explains the phase transitions in the KM (cf.~\cite{ChiMed17b}).
Importantly, we can trace the role of the network topology in the loss of stability of the incoherent state 
and  emerging spatial patterns.
In the remainder of this section, we informally review some of the results of \cite{ChiMed17a, ChiMed17b},
which will be used below.

The key ingredient in the analysis of the KM on graphs, which was not used in the analysis of the original KM,
is the 
the following kernel operator 
$\mathbf{W}:L^2(I)\to L^2(I)$ 
\be\lbl{def-W}
\bW[f](x)=\int_\R W(x,y) f(y) dy.
\ee
Recall that $W\in L^2(I^2)$ is a symmetric function. Thus, $\bW:~L^2(I^2)\to L^2(I^2)$ is a self-adjoint
compact operator. The eigenvalues of $\bW$, on the one hand, carry the information about the structure 
of the graphs in the sequence $\{\Gamma_n\}$, on the other hand, they appear in the stability analysis
of the incoherent state. Thus, through the eigenvalues of $\bW$ we can trace the relation between
the structure of the network and the onset of synchronization in the KM on graphs.

Since $\bW$ is self-adjoint and compact, the eigenvalues of $\bW$ are real and bounded with the only
possible accumulation point at $0$. 
In all examples considered in this paper, the largest eigenvalue $\mu_{max}>0$ and the 
smallest $\mu_{min}\le 0$. The linear stability analysis in \cite{ChiMed17a} shows that the incoherent state
is stable for $K\in [K_c^-, K_c^+],$ where 
$$
K_c^-= {2\over \pi g(0) \mu_{min}}\quad\mbox{and}\quad K_c^+= {2\over \pi g(0) \mu_{max}}.
$$

Except for small-world graphs, for all other graphs considered below, the smallest eigenvalue
$\mu_{min}=0$.  Thus, the incoherent state in the KM on these graphs is stable for $K\le K_c^+$.
In this section, we comment on the bifurcation at $K_c^+$ and postpone the discussion of the 
bifurcation at $K_c^-$ until Section~\ref{sec.SW}, where we deal with the KM on small-world graphs.

The analysis of the bifurcations in the KM on graphs
relies on the appropriate generalization of the Kuramoto's 
order parameter (cf.~\cite{ChiMed17a})
\be\lbl{dWorder}
h_n(t)=(h_{n1}(t), h_{n2}(t), \dots, h_{nn}(t)),\quad h_{ni}(t)=\frac{1}{n}\sum^n_{j=1}W_{n,ij}e^{i \theta _{nj}(t)}, 
\; i\in [n],
\ee
and its continuous counterpart
\be\lbl{Worder}
h(t,x) =\int_I\int_\R\int_\T W(x,y) e^{\1 \theta} \rho(t,\theta,\omega,y) g(\omega) 
d\theta d\omega dy.
\ee
The value of the continuous order parameter \eqref{Worder} $h(t,x)$ carries the information
about the (local) degree of coherence at point $x$ and time $t$. It is adapted to a given network 
connectivity through the kernel $W$.

The main challenge in the stability analysis of the incoherent state lies in the fact that
for $K\in [K_c^-, K_c^+]$ the linearized operator has continuous spectrum on the imaginary 
axis and no eigenvalues (cf.~\cite{ChiMed17a}). To overcome this difficulty, in \cite{ChiMed17b}, the generalized 
spectral theory was used to identify generalized eigenvalues responsible for the instability
at $K_c^+$, and to construct the finite-dimensional center manifold. The explanation of these
results in beyond the scope of this paper. An interested reader is referred to \cite{ChiMed17b}
for details. Here, we only present the main outcome of the bifurcation analysis. The center manifold
reduction for the order parameter near $K_c^+$ yields a stable branch of solutions:
\be\lbl{fork}
h_\infty (x, K)=
{ g(0)^2 \pi^{3/2} \over \sqrt{-g^{''}(0) }} \mu_{max}^{3/2} 
C(x) \sqrt{K-K_c^+} +o(\sqrt{K-K_c^+}), \quad 0<K-K_c^+\ll 1,
\ee 
where $C(x)$ is defined through the eigenfunction of $\bW$ corresponding to 
$\mu_{max}$ (see \cite{ChiMed17b} for details).
Formula \eqref{fork} generalizes the classical Kuramoto's formula describing
the pitchfork bifurcation in the all-to-all coupled model to the KM on graphs.
The network structure enters into the description of the pitchfork bifurcation
through the largest eigenvalue $\mu_{max}$ and the corresponding eigenspace.

\section{Results}\lbl{sec.results}
\setcounter{equation}{0}

In this section, we present numerical results elucidating some of the implications of the 
bifurcation analysis outlined in the previous section.

\subsection{Graphs approximating the complete graph}

\begin{figure}
\begin{center}
\textbf{a} \includegraphics[height=1.6in,width=2.0in]{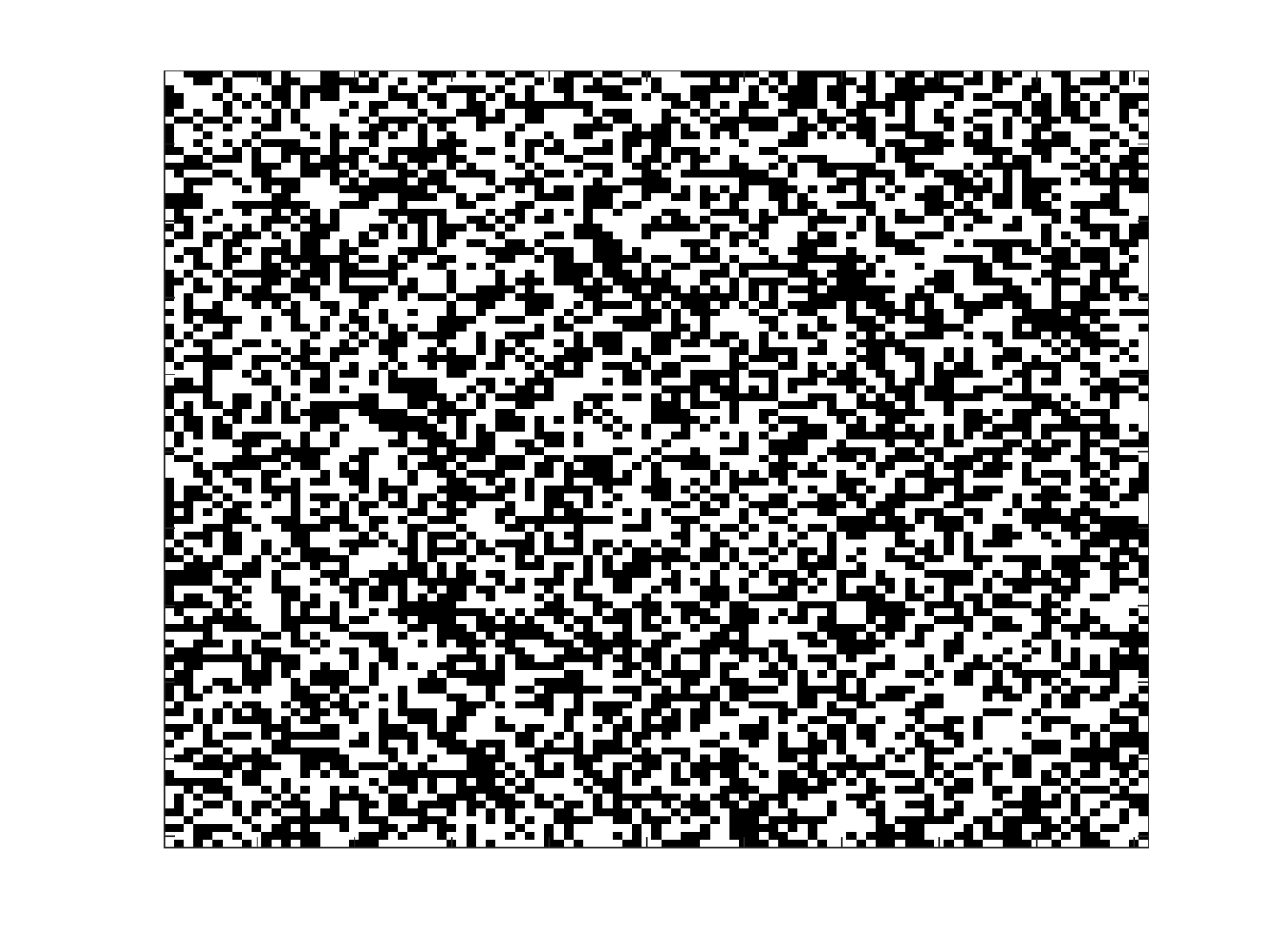}
\textbf{b} \includegraphics[height=1.6in,width=2.0in]{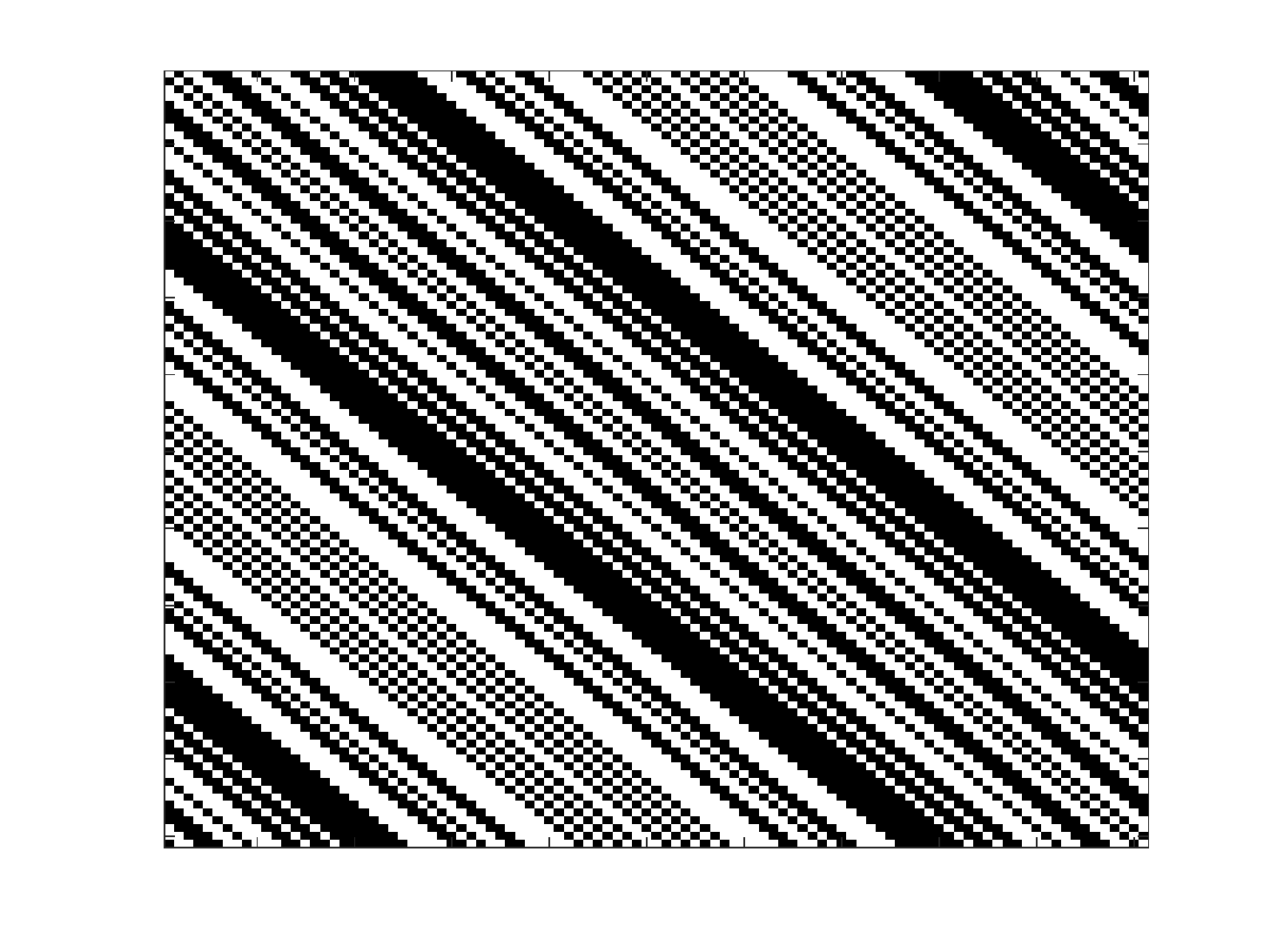}
\end{center}
\caption{Pixel pictures of Erd\H{o}s-R{\' e}nyi (\textbf{a}) and Paley (\textbf{b}) 
graphs.}
\label{f.ER-Paley}
\end{figure}

Our first set of examples deals with the KM on Erd\H{o}s-R{\' e}nyi and Paley  graphs (see Figure~\ref{f.ER-Paley}). 
The former is a random 
graph, whose edges are selected at random from the set of all pairs of vertices with fixed probability $p$
(see Example~\ref{ex.Wrandom}). Before giving
the definition of the Paley graphs, we first recall the definition of a Caley graph on an additive cyclic group
$\Z_n=\Z/n\Z$.
\begin{df}\lbl{df.Caley} Let $S$ be a symmetric subset of $\Z_n$ (i.e., $S=-S$). Then 
$\Gamma_n=\langle V(\Gamma_n), E(\Gamma_n)\rangle$ is a Caley graph if 
$V(\Gamma_n)=\Z_n$ and $\{a,b\}\in E(\Gamma_n),$ if $b-a\in S$. Caley graph
$\Gamma_n$ is denoted $\Cay(\Z_n,S)$.
\end{df}  
 
\begin{df}
Let  $n= 1\pmod 4$ be a prime and denote 
$$
Q_n=\{ x^2 \pmod{n}:\; x\in \Z_n/\{0\} \}.
$$  
$Q_n$ is viewed as a set (not multiset), i.e., each element has multiplicity $1$. Then $Q_n$ is 
a symmetric subset of $\Z_n^\times$ and $|Q_n|=2^{-1}(n-1)$ (cf. \cite[Lemma~7.22]{KreSha11}).
$P_n=\Cay(\Z_n, Q_n)$ is called a Paley graph \cite{KreSha11}.
\end{df}

\begin{figure}
\begin{center}
\textbf{a} \includegraphics[height=1.6in,width=2.0in]{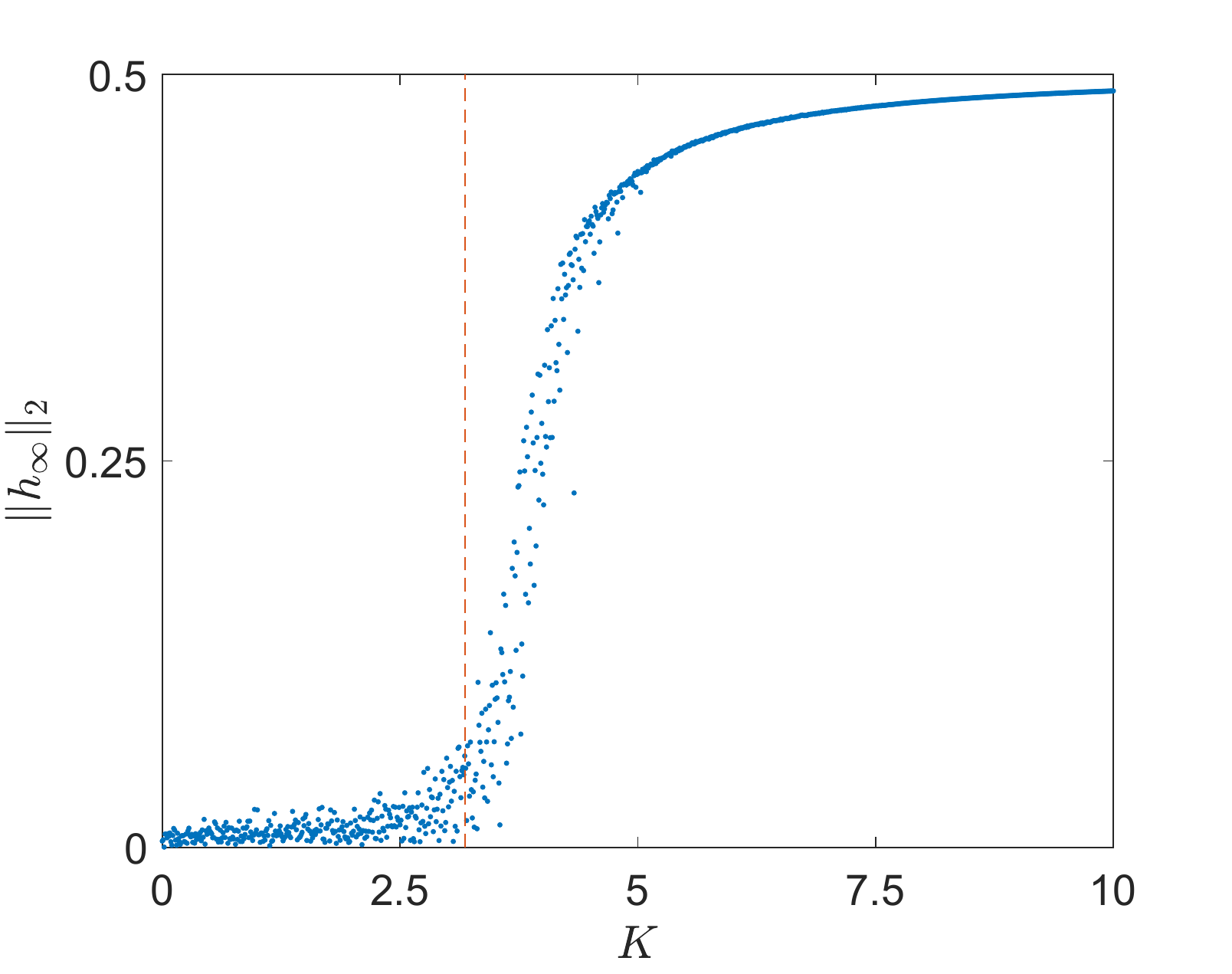}
\textbf{b} \includegraphics[height=1.6in,width=2.0in]{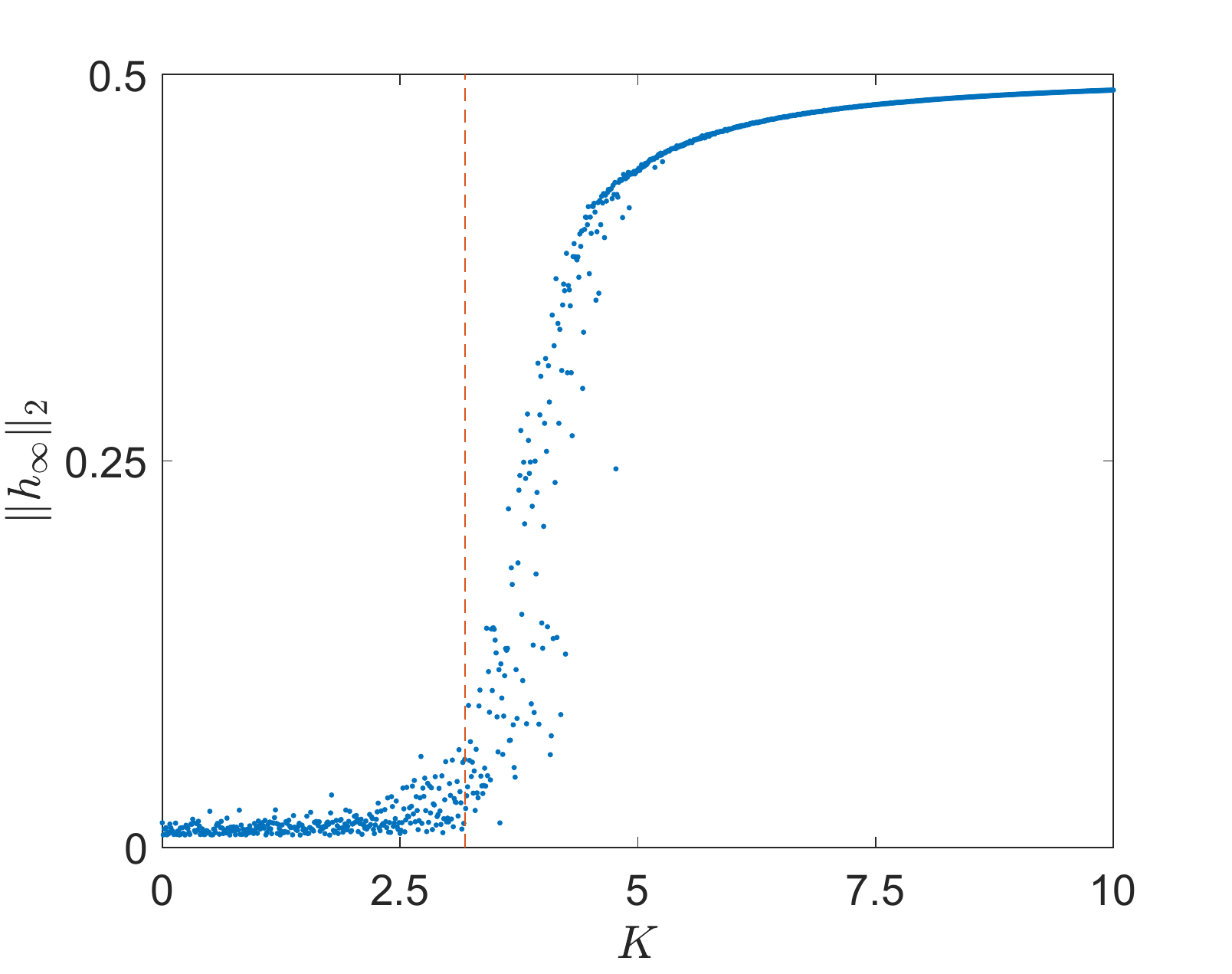}
\textbf{c} \includegraphics[height=1.6in,width=2.0in]{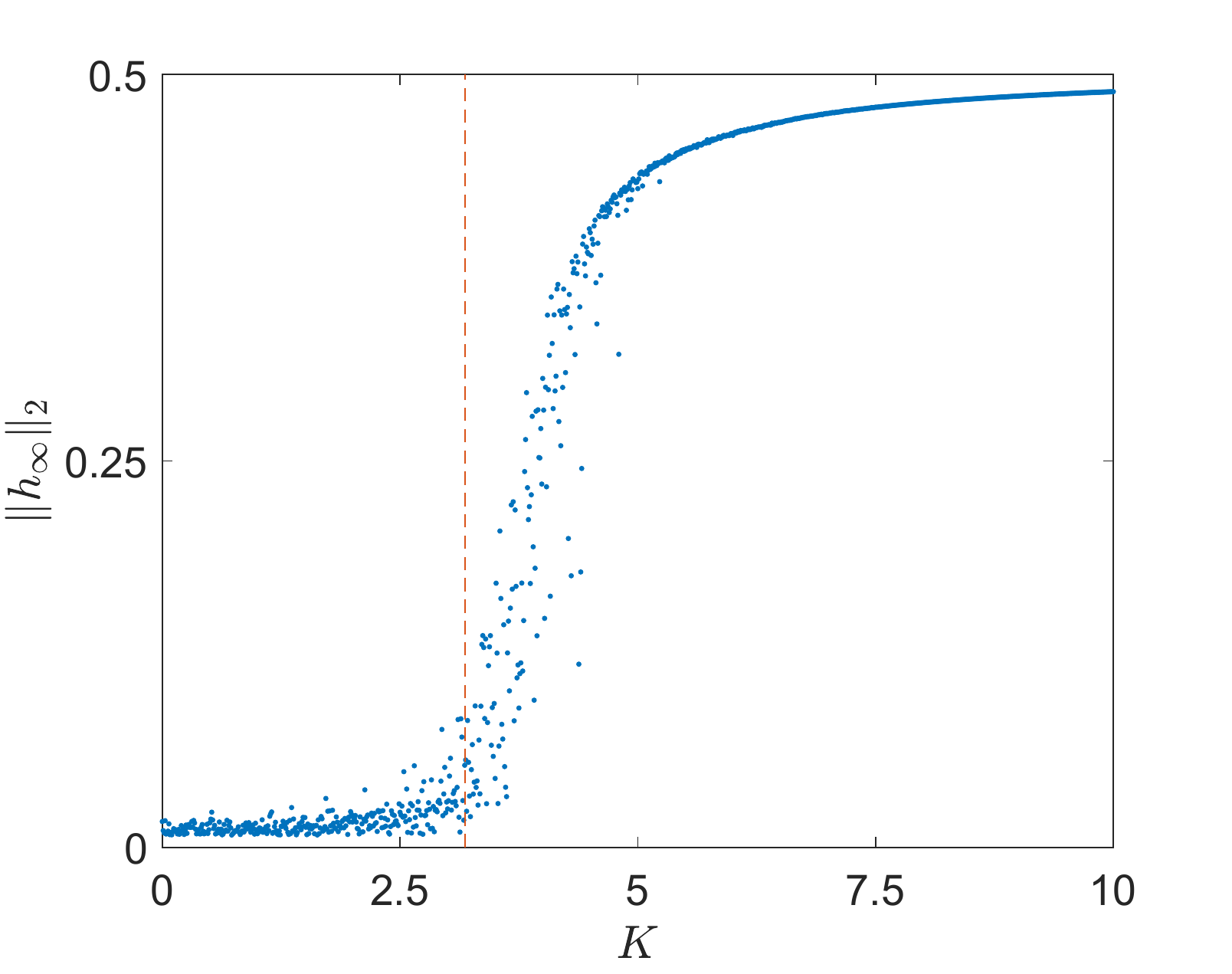}
\end{center}
\caption{A large time asymptotic value of the order parameter \eqref{dWorder} 
(in a scaled $\mathit{l}^2$-norm) of is plotted for the  KM on complete (\textbf{a}),
Erd\H{o}s-R{\' e}nyi (\textbf{b}), Paley (\textbf{c}). The transition from $0$ to $1/2$ takes
place in the same region of $K$ in all three plots. The critical value of the coupling strength
$K_c^+\approx 3.2$ is the same for all three models.}
\label{f.complete}
\end{figure}

As all Caley graphs, Paley graphs are highly symmetric (Figure~\ref{f.ER-Paley}\textbf{b}). 
However, Paley graphs also belong to
pseudorandom graphs, which share many asymptotic properties with Erd\H{o}s-R{\' e}nyi graphs with 
$p=1/2$ \cite{AlonSpencer}.
 In particular, both Erd\H{o}s-R{\' e}nyi and Paley graphs have constant graph limit $W_{1/2}\equiv 1/2$ 
as $n\to\infty$. The mean field limit for the KM for Erd\H{o}s-R{\' e}nyi and Paley graphs have the same 
kernel constant kernel $W=W_{1/2}\equiv 1/2$ (cf.~\eqref{MF}), exactly as in the case of the
 KM on weighted complete graphs. Since $\mu_{max}(W_{1/2})=1/2$, the transition to synchronization takes place
at
$$
K_c^+={4\over \pi g(0)}\approx 3.183
$$
for all three networks (see Figure~\ref{f.complete}). Here and in the numerical experiments below,
$g(x)=(2\pi)^{-1/2} e^{-x^2/2}$.

\subsection{Bipartite and disconnected graphs}
Our second set of examples  features another pair distinct networks, which are not pseudorandom, and yet  have
the same synchronization threshold.
\begin{figure}
\begin{center}
\textbf{a} \includegraphics[height=1.6in,width=2.0in]{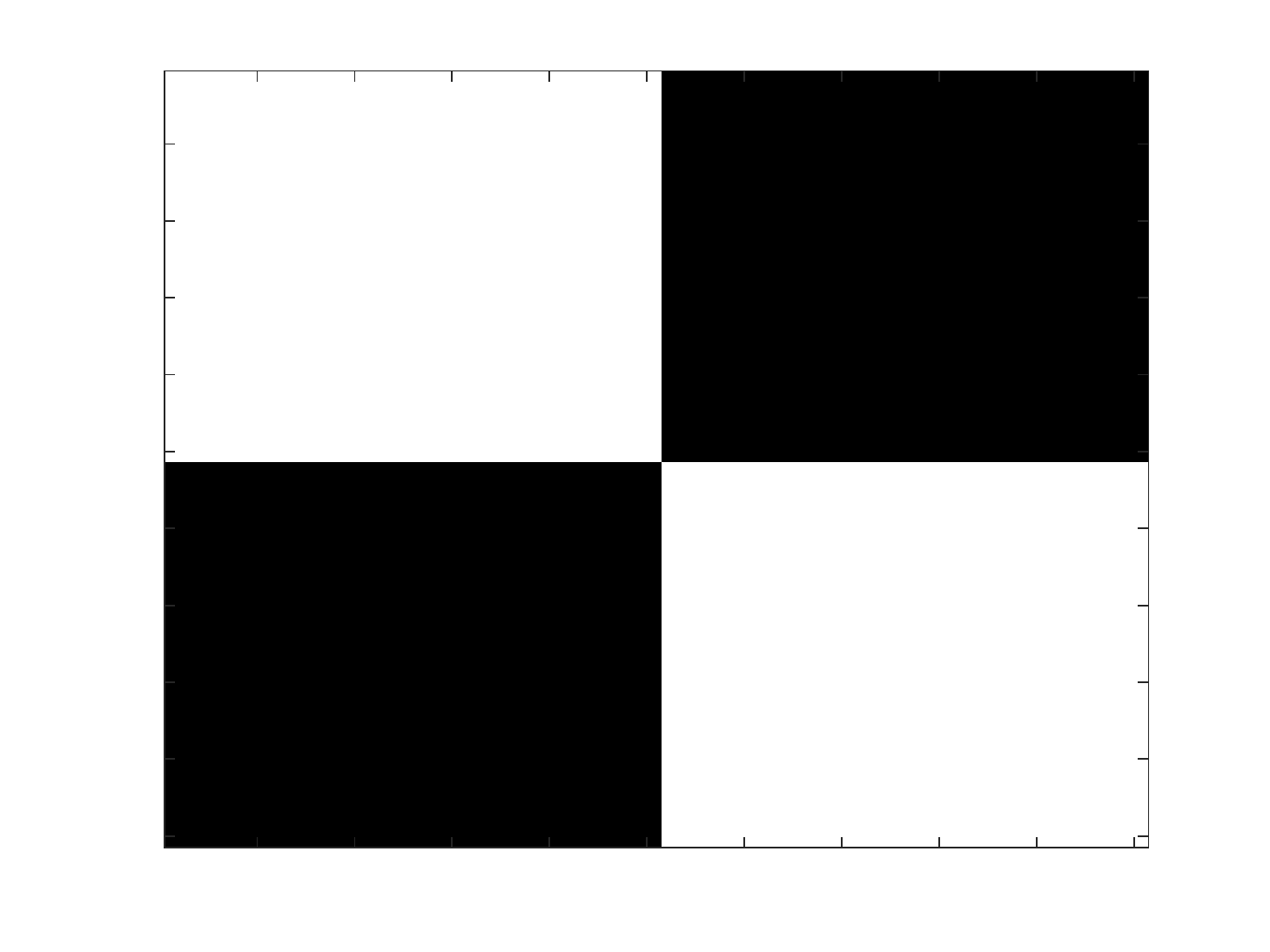}
\textbf{b} \includegraphics[height=1.6in,width=2.0in]{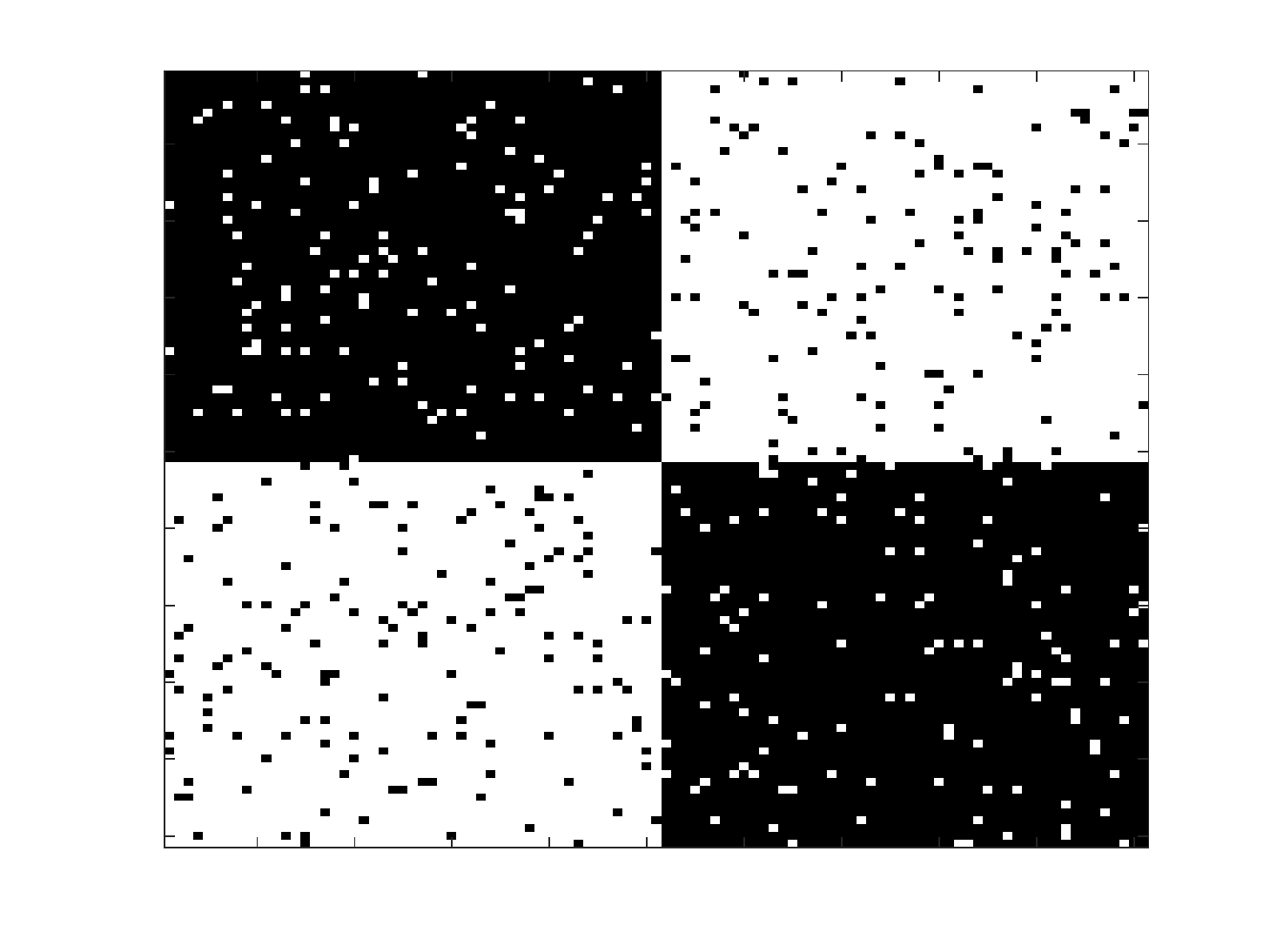}
\end{center}
\caption{Pixel pictures of the bipartite graph (\textbf{a}) and stochastic block graph with two weakly
connected components (\textbf{b}). 
graphs.}
\label{f.bi-disc}
\end{figure}

Here, we consider the KM on bipartite graphs $B_{n,n}=H(2n, W_{bp})$ (see Figure~\ref{f.bi-disc}\textbf{(a)})
\begin{equation*}
W_{bp}(x,y)=\left\{\begin{array}{ll}
1, & [0,0.5)\times [0.5, 1]\cup [0.5,1]\times [0,0.5),\\
0, & \mbox{otherwise.}
\end{array}\right.
\end{equation*}
Next, we consider the family of graphs interpolating between the disconnected graph with two equal
components and a weighted complete graphs, $\Gamma_{n,\alpha}=H_r(n,W_{d,\alpha})$:
\begin{equation*}
W_{d,\alpha}(x,y)=\left\{\begin{array}{ll}
1-\alpha, & [0,0.5)\times [0,0.5)\cup [0.5,1]\times [0.5, 1],\\
\alpha, & \mbox{otherwise.}
\end{array}\right.,\quad \alpha\in [0,0.5].
\end{equation*}

\begin{figure}
\begin{center}
\textbf{a} \includegraphics[height=1.6in,width=2.0in]{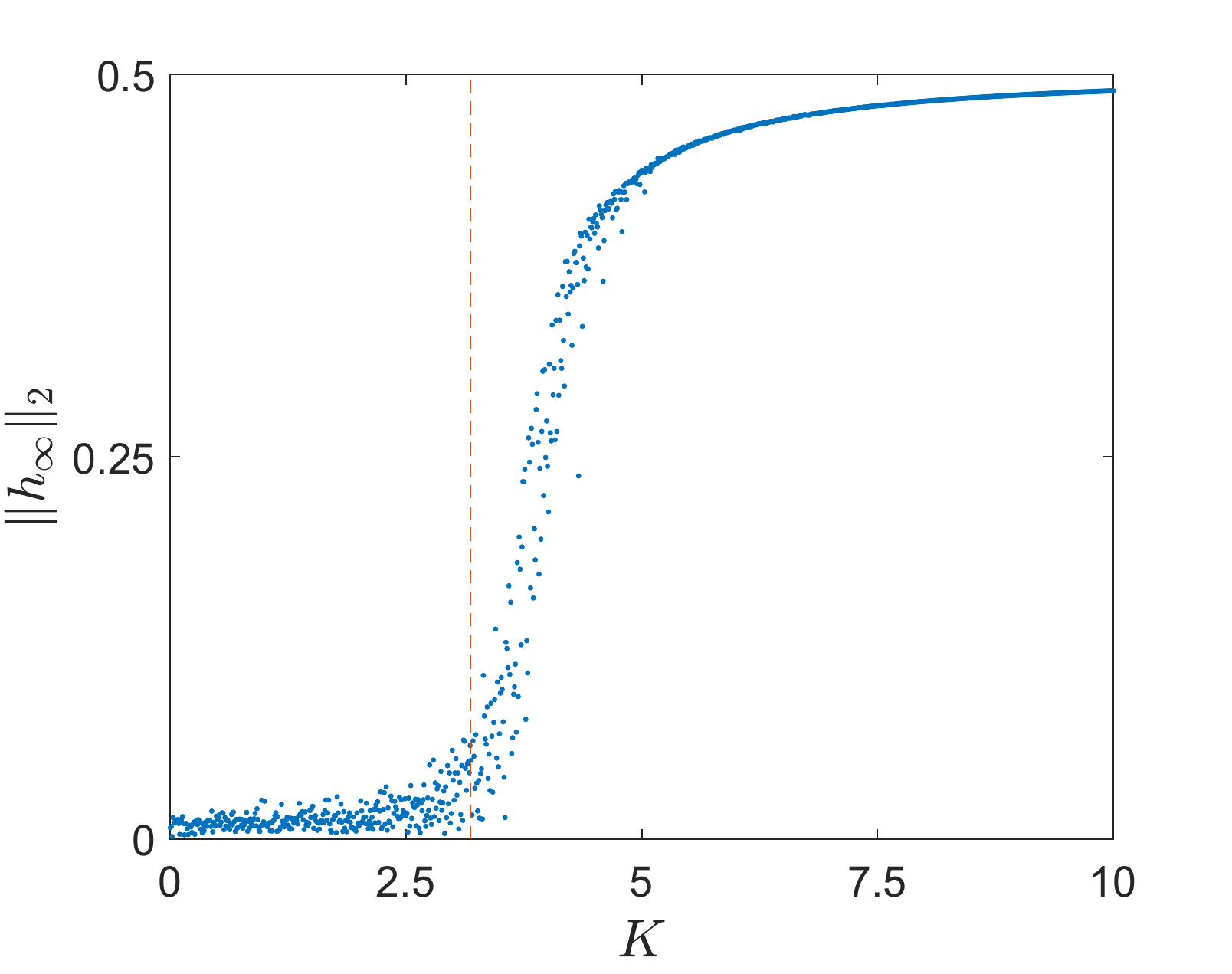}\quad
\textbf{b} \includegraphics[height=1.6in,width=2.0in]{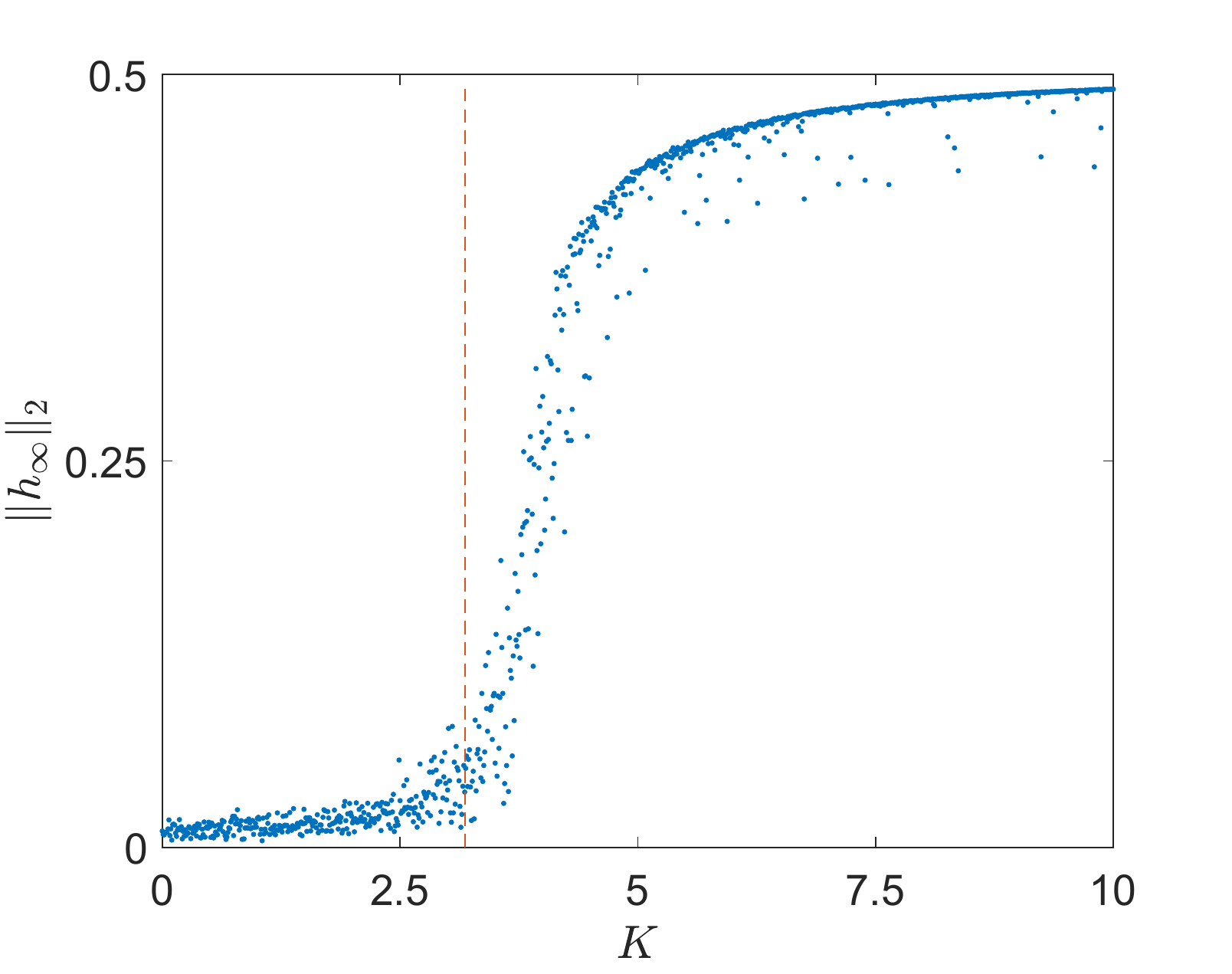}
\end{center}
\caption{The onset of synchronization in the KM on bipartite graph (\textbf{a}), 
and on the stochastic block graph with two weakly connected components (\textbf{b}) with $\alpha=.05$. 
}
\label{ff.bp}
\end{figure}

A simple calculation yields largest positive eigenvalues for each of these families 
of graphs
$$
\mu_{max}(W_{bp})=\mu_{max} (W_{d,\alpha})=0.5 \quad \forall \alpha\in [0, 0.5].
$$
Thus, as in the first set of examples the KM on $B_{n,n}$ and $\Gamma_{n,\alpha}$
undergoes transition to synchronization at $K_c^+=4\left(\pi g(0)\right)^{-1}$ 
(see Figure~\ref{ff.bp}).

\subsection{Power law graphs}
\begin{figure}
\begin{center}
\includegraphics[height=1.6in,width=2.0in]{powerlaw.pdf}
\includegraphics[height=1.6in,width=2.0in]{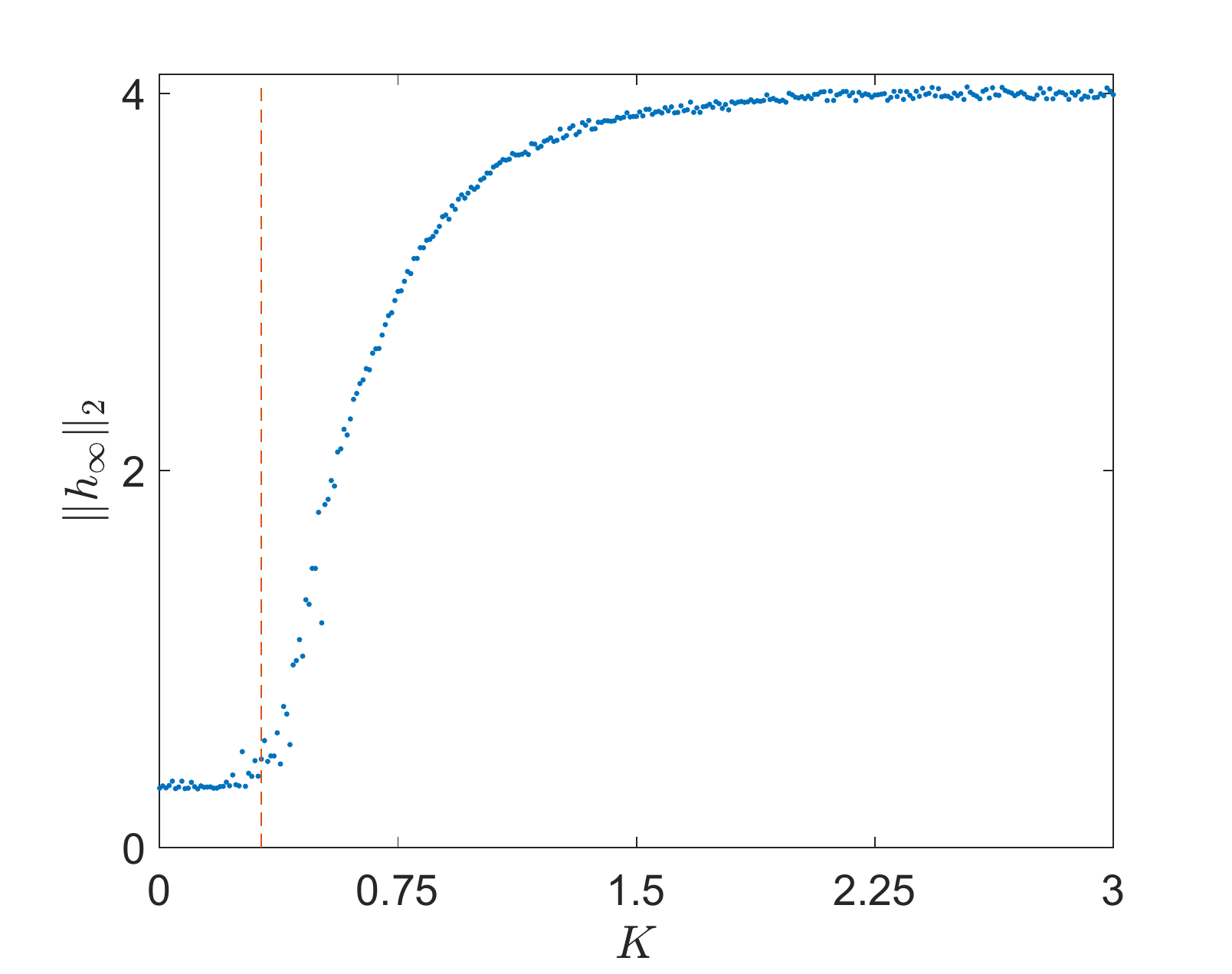}
\end{center}
\caption{The onset of synchronization in the KM on power-law graphs. The parameter values $\gamma=.4$
and $\beta= .6$ were used in these numerical experiments.
}
\label{f.pl}
\end{figure}

Next we turn to the KM on power law graphs. To generate graphs with power law degree
distribution, we use the method of sparse W-random graphs \cite{BCCZ}.
Specifically, $\Gamma_n=H_r(n,W^\gamma, n^{-\beta})$ with graphon $W^\gamma$ defined
in \eqref{Walpha}. For graphs constructed using this method, it is known that
the expected degree of node $i\in [n]$ is 
$O(n^{1+\gamma-\beta})$ and the edge density is $O(n^{-\beta})$
(cf.~\cite{KVMed17}). Thus, $\Gamma_n=H_r(n,W^\gamma, n^{-\beta})$ is a family of sparse
graphs with power law degree distribution. The mean field limit for the KM on $\Gamma_n$
is given by \eqref{MF} with $W:=W^\gamma$. The analysis of the integral operator 
$\bW$ (cf. \eqref{def-W}) with kernel $W:=W^\gamma$ shows that it has a single nonzero
eigenvalue (cf.~\cite{MedTan17})
$$
\mu_{max}=(1-2\gamma)^{-1}, \quad \gamma\in (0, 0.5).
$$
Thus, the synchronization threshold for the KM on power law graphs is
$$
K^+_c= {2(1-2\gamma)\over \pi g(0)}.
$$
Note that as  $\gamma\to 1/2,$ $K_c^+$ tends to $0$. Thus, the KM on power
law graphs features remarkable synchronizability despite sparse connectivity.
Numerics in Figure~\ref{f.pl} illustrate the onset of synchronization in power
law networks.

\section{Small-world graphs}\lbl{sec.SW}
\setcounter{equation}{0}
Small-world graphs are obtained by random rewiring of regular $k$-nearest-neighbor
 networks \cite{WatStr98}. 
Consider a graph on $n$ nodes arranged around a circle, with each node connected to $k=\lfloor rn\rfloor$
neighbors from each side for some $r\in (0,0.5)$. With probability
$p\in [0, 0.5]$ each of these edges connecting $k$ neighbors from each
side, is then replaced  by a random long-range (i.e., going outside the $k$-neighborhood) edge.
The pixel picture of a representative
small-world graph is shown in  Figure~\ref{f.graphs}\textbf{a}. 
A family of small-world graphs parametrized by $p$ interpolates between regular $k$-nearest-neighbor
graph ($p=0$) and a fully random Erd\H{o}s-R{\' e}nyi graph ($p=0.5$). 
For intermediate values of $p,$ small-world graphs combine features of regular and random connectivity,
just as seen in many real-world networks \cite{WatStr98}. 

It is convenient to interpret a small-world graphs as a $W$-random graph $S_{n,p,r}=H_r(n,W_{r,p})$ 
(cf.~\eqref{def-Wpr}).
The mean field limit of the KM on $S_{n,p,r}$ is then given by \eqref{MF} with $W:=W_{p,r}$.
The corresponding kernel operator is given by
\be\lbl{convolution}
\bW_{p,r}(f)= \int_{\T_1} K_{p,r}(\cdot-y)f(y)dy,\quad \T_1:=\R/\Z.
\ee
 where $K_{p,r}(x)=(1-p) \mathbf{1}_{|x|\le r}(x)+p \mathbf{1}_{|x|>r}(x).$

Using \eqref{convolution}, we recast the eigenvalue problem for 
$\textbf{W}_{p,r}$ as follows
\be\lbl{EVP}
K_{p,r} \ast v=\mu v.
\ee
By taking the Fourier transform of the both sides of \eqref{EVP}, we find 
the eigenvalues of $\textbf{W}_{p,r}$ 
\be\lbl{zeta-k}
\mu_k =\int_{\T}\! K_{p,r}(x)e^{-2\pi i kx}dx =\left\{ \begin{array}{ll} (2r+p-4rp), & k=0,\\
(\pi k)^{-1} (1-2p) \sin (2\pi kr), & k\in\Z/\{ 0\}.
\end{array}
\right.
\ee
The corresponding eigenvectors are $w_k= e^{\1 2\pi kx}, k\in \Z.$ 
Note that $\mu_k=\mu_{-k},$ since $K_{p,r}$ is even (cf.~\eqref{zeta-k}). 
Thus, the eigenspace corresponding to $\mu=\mu_0$ is spanned by $w_0=1$.
For $\mu\neq \mu_0,$ 
the eigenspace corresponding to $\mu$ is spanned by
\be\lbl{2d-EV}
w_k=e^{\1 2\pi kx}\quad \mbox{and}\quad w_{-k}=e^{-\1 2\pi kx}, \quad k\in I_\mu,
\ee
where $I_\mu=\{k\in \N:~\mu_k=\mu\}.$
The largest positive eigenvalue of $\textbf{W}_{p,r}$ is
$$
\mu_0= \int_I K_{p,r} (x) dx =2r+p-4pr.
$$
Therefore, the onset of synchronization in the KM on small-world graphs takes place at
$$
K_c^+={2\over \pi g(0) (2r+p-4pr)}
$$
(see Figure~\ref{f.2}\textbf{a}).

Importantly, $\textbf{W}_{p,r}$ also has negative eigenvalues. Since $\bW_{p,r}$ is a compact
operator, $0$ is the only accumulation point of the spectrum of $\bW_{p,r}$.
Thus, there is the smallest (negative) eigenvalue
$\mu_{min}=\min\{\mu_k:\; k\in \N\}$.
Let $q\in \N$ be such that $\mu_q=\mu_{min}$. Assuming that the multiplicity of $\mu_{min}$
is $2$, the center manifold reduction performed for the order parameter in \cite{ChiMed17b}
yields the following stable branch bifurcating from the incoherent state ($h\equiv 0$) 
at $K=K^-_c$
\be\lbl{hinfty}
h_\infty (x,K)= \sqrt{\kappa\over \beta} e^{\pm \1 2\pi q (x+\phi)} + o(\sqrt{\kappa}),\quad 
0<\kappa:=K_c^--K\ll 1,\; \beta:=-\frac{\pi g''(0) (K_c^-)^4\mu_{min}}{16},
\ee
where $\phi\in [0,1)$ depends on the initial data.

Equation \eqref{hinfty} implies that at $K=K_c^-$ the KM on small-world graphs undergoes a pitchfork bifurcation.
Unlike the bifurcation at $K=K^+_c$ considered earlier, in the present case the center manifold is two-dimensional.
It is spanned by $w_{\pm q}=e^{\pm\1 2\pi qx}$. In the remainder of this section, we analyze stable spatial patterns bifurcating
from the incoherent state at $K=K_c^-$.
To this end, we rewrite the velocity field \eqref{def-V}, using the order parameter:
\be\lbl{V-h}
V(t,u,\omega,x) = \omega +{K\over 2\1}\left( e^{-\1 u} h(t,x)- e^{\1 u} \overline{h(t,x)}\right).
\ee
The velocity field in the stationary regime is then given by
\be\lbl{Vinfty}
V_\infty (u,\omega,x) =\omega +{K\over 2\1} \left(  e^{-\1 u} h_\infty(x,K)- e^{\1 u} \overline{h_\infty(x,K)}\right).
\ee
Using the polar form of the order parameter
\be\lbl{polar}
h_\infty(x,K_c^--\kappa)=R_\infty(x,K) e^{i\Phi(x)}, \quad R_\infty(x,K)=\sqrt{\kappa\over \beta}+o(\sqrt{\kappa}),
\Phi(x)=\pm 2\pi q (x+\phi),
\ee
from \eqref{Vinfty} we obtain
\be\lbl{Vinfty-a}
V_\infty (u,\omega,x)= \omega - (K_c^--\kappa) \sqrt{\kappa\over \beta} \sin\left(u- \Phi(x) \right)+ 
o(\sqrt{\kappa}).
\ee

\begin{figure}
	\begin{center} 
\textbf{a}\includegraphics[height=1.6in,width=2.0in]{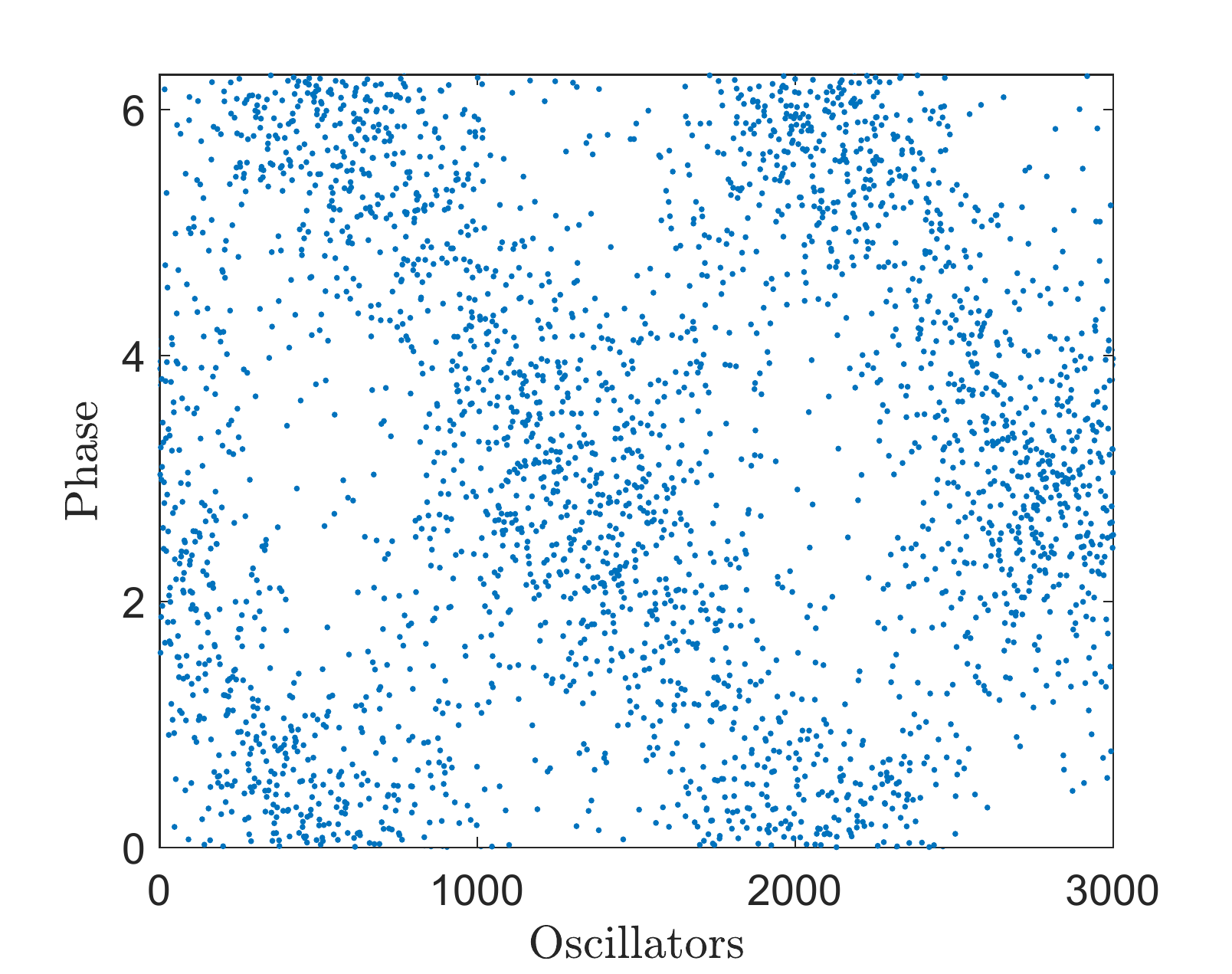}
\textbf{b}\includegraphics[height=1.6in,width=2.0in]{k_-36_r_pt3.pdf}
\textbf{c}\includegraphics[height=1.6in,width=2.0in]{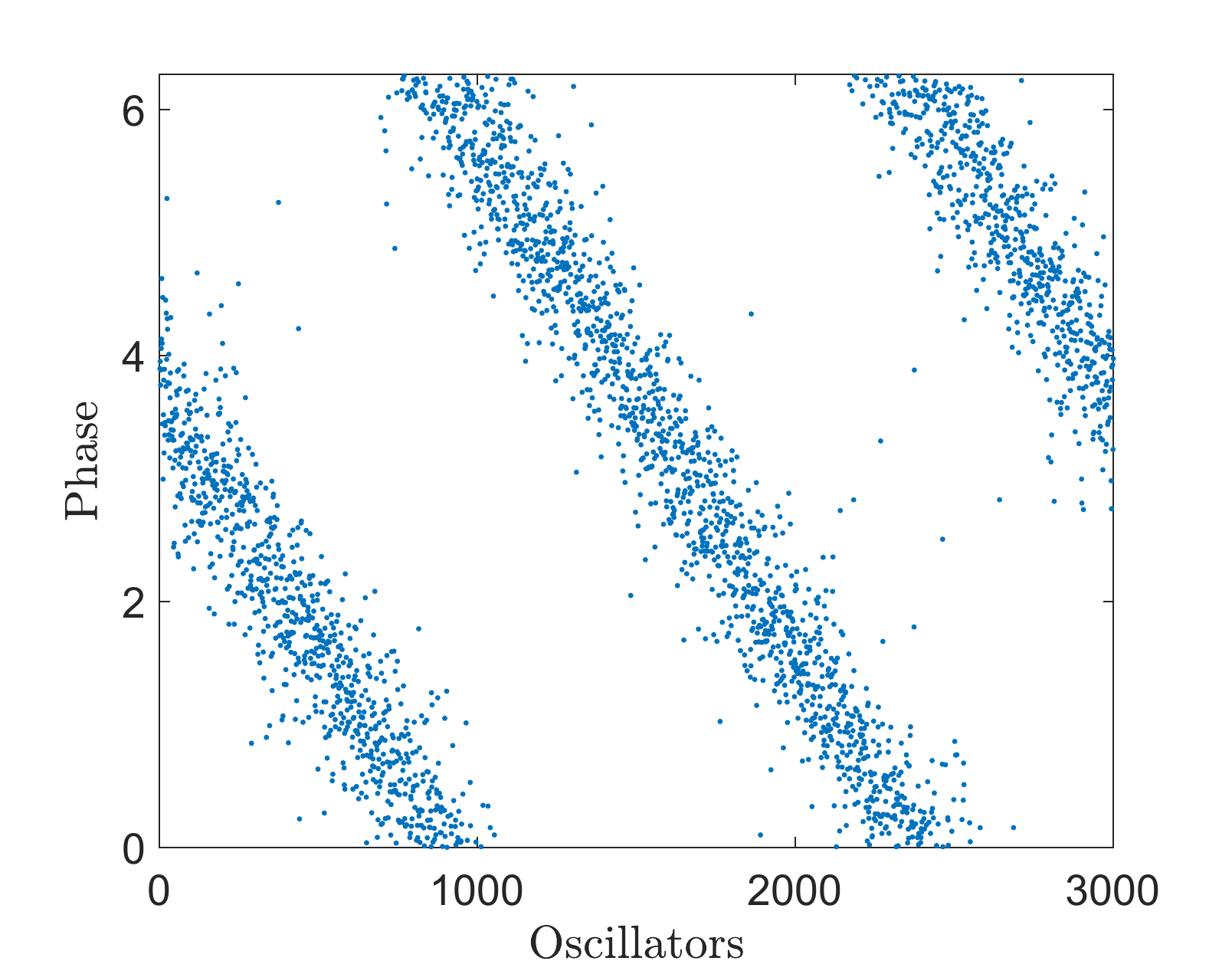}
\end{center}
\caption{The formation of 2-twisted states for decreasing values of coupling strength $K$. 
From {\bf a} to {\bf c}, we have, respectively, $K=-32$, $K=-36$, and $K=-50$. In all simulations $p=.2$ and $r=.3$.}
\label{f.twist}
\end{figure}

\begin{figure}
\begin{center}
	\includegraphics[height=1.6in,width=2.0in]{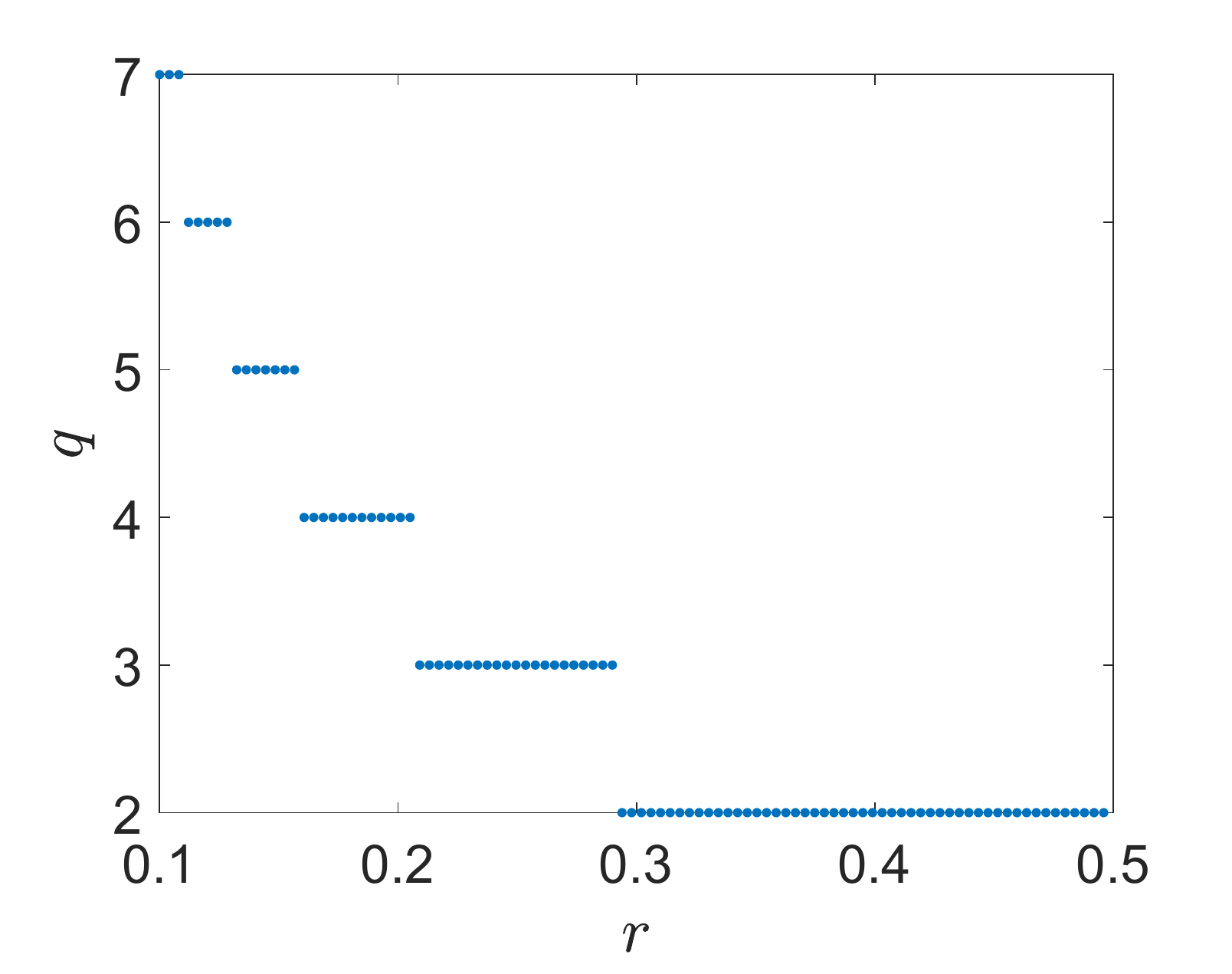}
\caption{The winding number $q$ of the twisted states bifurcating from
  the incoherent state as a function of the connectivity range $r$.
}
\label{f.winding}
\end{center}
\end{figure}

Since $\rho_\infty$ is a steady state solution of the \eqref{MF}, it satisfies 
\be \lbl{steadyV}
\p_u\left\{  V_\infty \rho_\infty \right\}=0.
\ee
From \eqref{Vinfty-a} and \eqref{steadyV}, we have
\be\lbl{partially-locked}
\rho_\infty (u ,\omega ,x) =
\left\{  \begin{array}{ll} \delta \left(u -\Phi(x) -\arcsin \left(\omega {\sqrt{ \beta}\over 
K\sqrt{\kappa}}\right)
+o(\sqrt{\kappa})\right), & |\omega |\le KR_\infty(x,K),\\
 \frac{1}{2\pi} \frac{\sqrt{\omega ^2 - K^2 R_\infty(x,K)^2}} {| \omega - KR_\infty(x,K)
\sin (u -\Phi (x))|}, &|\omega |>KR_\infty(x,K).
\end{array}\right.
\ee
Here, $K=K_c^--\kappa$ and $\delta$ stands for the Dirac delta function. 
Equation \eqref{partially-locked} describes partially phase locked solutions. They combine phase locked 
oscillators (type I)
\be\lbl{twist}
u = \pm 2\pi q(x+\phi) + Y(\omega, \kappa)+ O\left(\sqrt{\kappa}\right), \quad 
Y(\omega, \kappa)=\arcsin \left(\omega {\sqrt{ \beta}\over K\sqrt{\kappa}}\right)
\ee
and those whose distribution is given by the second line in \eqref{partially-locked} (type II).
The oscillators of the first type form $q$-twisted states subject to noise $Y$. Note that 
$Y$ is a function of $\omega$ and, therefore, is random.  
Just near the bifurcation, where  $0<\kappa\ll 1,$ the probability of $|\omega |>KR_\infty(x,K)$
is high and, thus, most oscillators are of the second type
(Figure~\ref{f.twist}\textbf{a}). However, as we move away from the bifurcations, we see the number
of the oscillators of the first type increase and the noisy twisted states become progressively more
distinct (Figure~\ref{f.twist} \textbf{b},\textbf{c}).
Thus, the bifurcation analysis of the KM  on small-world graphs 
identifies a family of stable twisted states (cf.~\cite{MedTan15b, WilStr06}). By changing parameters 
of small-world connectivity, one can control the winding number of the emerging twisted states
(see Figure~\ref{f.winding}).

\section{Discussion}\lbl{sec.discuss}
\setcounter{equation}{0}

In this paper, we selected several representative families of graphs to illustrate the link between
network topology and synchronization and pattern formation in the KM of coupled phase oscillators.
In particular, we showed that the transition to synchronization in the KM on pseudorandom graphs
(e.g., Erd\H{o}s-R{\' e}nyi, Paley, complete graphs) starts at the same critical value of the coupling
strength and proceeds in practically the same way.  The bifurcation plots shown in 
Figure~\ref{f.complete} are very similar, although plots for the KM on  Erd\H{o}s-R{\' e}nyi and Paley
graphs show more variability. The almost identical bifurcation plots for these models are due to the 
fact that all three models result in the same mean field equation. This means that the empirical 
measures generated by the trajectories of these models are asymptotically have the same limit as
$n\to\infty$. The differences seen in plots \textbf{a}-\textbf{c} of Figure~\ref{f.complete} are
due to the higher order moments, which are not captured by the mean field limit.
Other families of graphs considered in this paper include the bipartite graph and a family of 
stochastic block graphs interpolating between a graph with two disconnected components and   
Erd\H{o}s-R{\' e}nyi graph. The KM on all these graphs (including the disconnected graph) feature
the same transition to synchronization. Finally, we studied the bifurcations in the KM on small-world
and power law graphs due to their importance in applications. A remarkable feature of the KM on 
small-world graphs is the presence of the bifurcation leading to stable noisy twisted states. 
Twisted states were known as attractors in repulsively coupled KM with identical intrinsic 
frequencies \cite{WilStr06, Med14c, MedTan15b, MW17}. In this paper, we show that twisted
states (albeit noisy) can also be attractors in the KM with random intrinsic frequencies.
Furthermore, they bifurcate from the incoherent, i.e. fully random, steady state.

\begin{figure}
	\begin{center}
		\includegraphics[width = .3\textwidth]{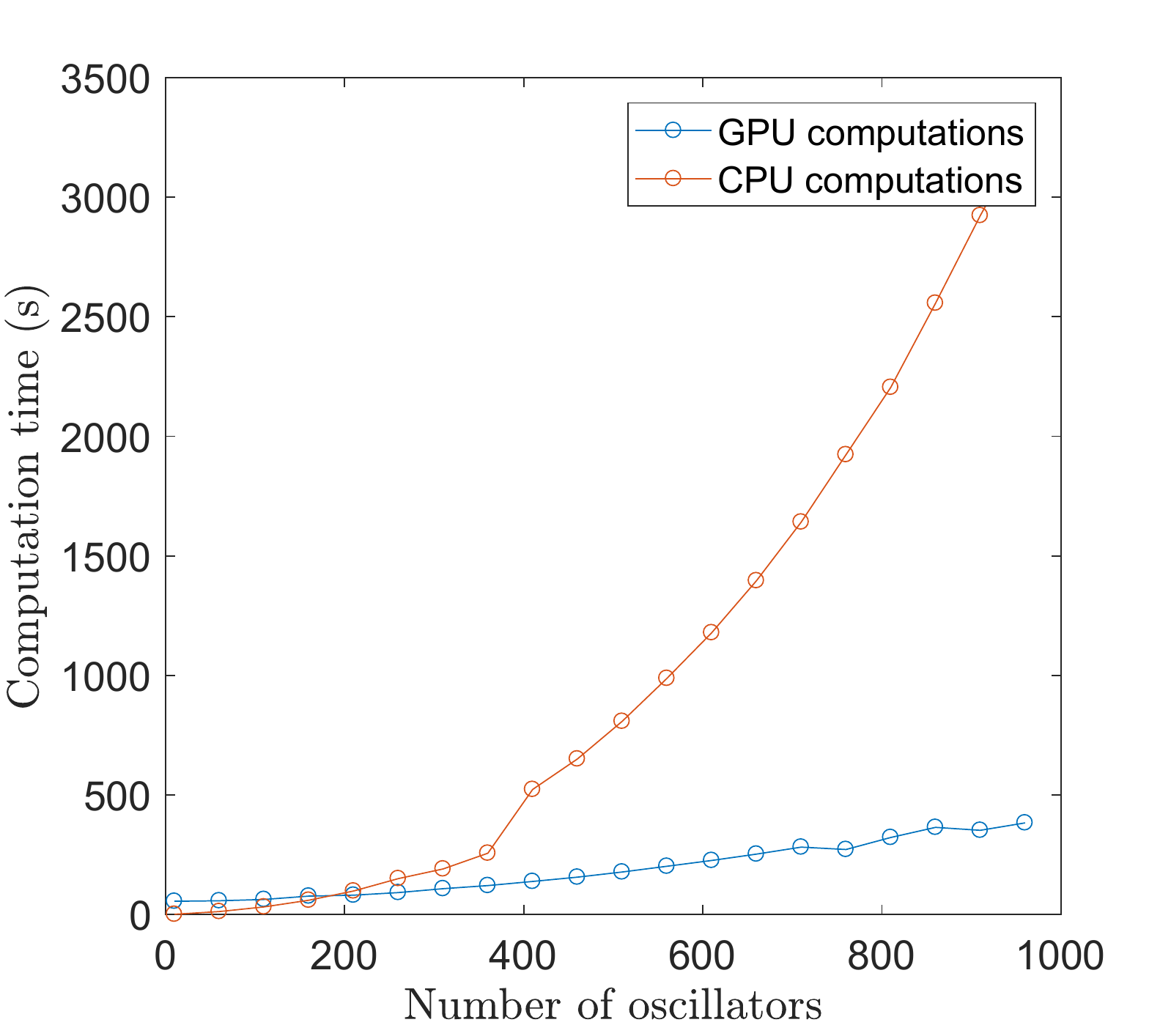}
	\end{center}
	\caption{Comparison of GPU and CPU computation time on Intel(R) Core i7-7700 CPU
 with NVIDIA Quadro K1200 GPU showing dramatic speed-up due to parallelization of computations.}
	\label{cpu_gpu}
\end{figure}

To perform numerical experiments presented in this paper we had to overcome several
challenges. Verification of the bifurcation scenarios in the KM required a large number
of simulations of large systems of ordinary differential equations with random coefficients.
Thus, the speed of computations was critical in this project.
All computations were completed in MATLAB utilizing GPU computations for dramatic speed-up 
in the computational time (compared to CPU computations), see Figure \ref{cpu_gpu}. 
Time steps were performed using Heun's method with $\Delta t = .01$, which was sufficient for 
stable simulation results.

Each simulation was initialized with a random state vector $(u_{n,1},\dots, u_{n,n})^T$ 
with each component independently chosen from a uniform distribution on $[0,2\pi)$ 
representing the phase of each oscillator. We additionally initialize a random vector of
 internal frequencies chosen from a normal distribution. Unless otherwise noted, in all 
simulations we take $n=4001$, which is a prime $1$ modulo $4$, so that the theory developed 
for Paley graphs applies. We observed that taking $T_{final}=20$ was sufficient for systems to 
exhibit synchronization or $q$-twisted states. 

Computationally solving \eqref{zeta-k} for the minimal negative eigenvalue is accomplished by observing the 
trivial bound $\mu_k>\frac{-1}{k\pi}$ and iteratively computing $\mu_{min}=\mu_{k_0}$ for $k_0\in \{1,\dots, M\}$
 for increasing values of $M$ until $\mu_{min}\leq \frac{-1}{M\pi}$.

\subsection*{Acknowledgments}

The work of the second author was supported in part by the NSF DMS 1715161.
The authors would like to thank Nicholas Battista for helpful discussions related to 
implementation of the numerical simulations.

\vfill\newpage
\bibliographystyle{amsplain}

\def\cprime{$'$} \def\cprime{$'$} \def\cprime{$'$}
\providecommand{\bysame}{\leavevmode\hbox to3em{\hrulefill}\thinspace}
\providecommand{\MR}{\relax\ifhmode\unskip\space\fi MR }
\providecommand{\MRhref}[2]{%
  \href{http://www.ams.org/mathscinet-getitem?mr=#1}{#2}
}
\providecommand{\href}[2]{#2}

\end{document}